\documentclass{article}
\usepackage[top=1in, bottom=1in, left=1in, right=1in]{geometry}
\usepackage{amsmath, amsfonts, amssymb, amsthm}
\usepackage{hyperref}
\hypersetup{
    pdfauthor={Timothy Eller},
    pdftitle={Chiral vector bundles},
    colorlinks=true,
    linkcolor=brown,
    citecolor=purple,
    urlcolor=blue
}
\usepackage{tikz}
\usetikzlibrary{arrows,chains,matrix,positioning,scopes}

\theoremstyle{plain}
\newtheorem{thm}{Theorem}[section]
\newtheorem{lem}[thm]{Lemma}
\newtheorem{cor}[thm]{Corollary}
\newtheorem{prop}[thm]{Proposition}
\newtheorem{conj}[thm]{Conjecture}

\theoremstyle{definition}
\newtheorem{rmk}[thm]{Remark}
\newtheorem{defn}[thm]{Definition}
\newtheorem{example}{Example}
\newtheorem*{note}{Notation}

\newcommand{\sembrack}[1]{[\![#1]\!]}
\newcommand{\R}{\textbf{R}}
\newcommand{\C}{\textbf{C}}
\newcommand{\Z}{\textbf{Z}}
\newcommand{\1}{\textbf{1}}
\newcommand{\id}{\operatorname{id}}
\newcommand{\Aut}{\operatorname{Aut}}
\newcommand{\Int}{\operatorname{Int}}
\newcommand{\ob}{\operatorname{ob}}
\newcommand{\End}{\operatorname{End}}
\newcommand{\res}{\operatorname{res}}
\newcommand{\Res}{\operatorname{Res}}
\newcommand{\supp}{\operatorname{supp}}
\newcommand{\spanof}{\operatorname{span}}

\title          {Chiral vector bundles}
\author         {Timothy Eller}

\begin {document}
\maketitle

\abstract
Given a smooth $G$-vector bundle $E\to M$ with a connection $\nabla$, we propose the construction of a sheaf of vertex algebras $\mathcal{E}^{ch(E,\nabla)}$, which we call a \textit{chiral vector bundle}. $\mathcal{E}^{ch(E,\nabla)}$ contains as subsheaves the sheaf of superalgebras $\Omega \otimes \Gamma (SE \otimes \Lambda E)$ and the sheaf of Lie algebras generated by certain endomorphisms of these superalgebras: $\nabla$, the infinitesimal gauge transformations of $E$, and the contraction operators $\iota_X$ on differential forms $\Omega$. Another subsheaf of primary importance is the chiral vector bundle $\mathcal{E}^{ch(M\times \C,d)}$, which is closely related to the chiral de Rham sheaf of Malikov et alii.

\tableofcontents

\section{Introduction}

\subsection{Overview}

The advent of the chiral de Rham complex \cite{MSV-1999} has introduced vertex algebras to the context differential and algebraic geometry. Namely, the vertex algebra analogs of the Heisenberg and Clifford algebras (known in other literature either as the $\beta\gamma$ and $bc$ systems, or fermionic and bosonic ghost systems) patch together to form a sheaf, over a smooth manifold $M$, that contains the classical de Rham complex as a subsheaf. In the past decade, several aspects of geometry have subsequently been lifted to their vertex algebra analogs, a process we might call \textit{chiralization}. Recent such examples include the chiralization of \textit{equivariant} cohomology by \cite{LL-2007} and its implications, such as \cite{tan-2009-03}; and the chiralization of differential operators on $M$ \cite{kapranov-2006}, among others.

In this paper, we continue this trend, extending \textit{chiralization} to vector bundles and their geometry.

\subsubsection{Relation to the chiral de Rham sheaf}
\label{sec:CDRrelationship}

Just as the ordinary de Rham complex resides inside the chiral de Rham sheaf as a subsheaf, our final goal is to take a $G$-vector bundle $E$ with connection $\nabla$ and exhibit the sheaf of sections of $\Omega \otimes \Gamma (SE \otimes \Lambda E)$ ($\Omega$ being the exterior algebra of differential forms, and $SE$ and $\Lambda E$ being the full symmetric and antisymmetric tensor algebras of $E$) as a subsheaf of $\mathcal{E}^{ch(E,\nabla)}$, a sheaf of vertex algebras which we call a \textit{chiral vector bundle}. In fact $\mathcal{E}^{ch(E,\nabla)}$ will contain the Lie algebra of endomorphisms of $\Omega \otimes \Gamma (SE \otimes \Lambda E)$ generated by $\nabla$, the infinitesimal gauge Lie algebra, and the odd Lie algebra of contraction operators $\iota_X$ ($X \in \mathfrak{X}$, the set of smooth vector fields), as well.

The construction of the vertex algebra $\mathcal{E}^{ch(E,\nabla)}(U)$, $U\subset M$, is outlined as follows. Restricting all objects to $U$, we enlarge the space $\Omega \otimes \Gamma E$ to the supersymmetric algebra $\mathfrak{a} = \Omega \otimes \Gamma (SE \otimes \Lambda E)$. We then capture the connection $\nabla$ in a robust Lie superalgebra $\mathfrak{s} \mathfrak{X}_{\nabla}$ of endomorphisms of the aforementioned superalgebra $\mathfrak{a}$.

This pair is combined into a single entity $\mathfrak{s} \equiv \mathfrak{s} (\mathfrak{a} \otimes \mathfrak{s} \mathfrak{X}_{\nabla}, \mathfrak{a})$ (the tensor products taken over the ring $\Omega^0$ of smooth functions) which we generically call a \textit{souped-up Lie algebra}. It is the semidirect sum Lie algebra of $\mathfrak{a} \otimes \mathfrak{s} \mathfrak{X}_{\nabla}$ and $\mathfrak{a}$, also with the structure of a left module for $\mathfrak{a}$. The souped-up Lie algebra $\mathfrak{s}$ is a classical construction; we have not mentioned vertex algebras so far.

The impetus for constructing the souped-up Lie algebra is that a vertex algebra generated by the underlying Lie algebra \textit{almost} has the structure of a souped-up Lie algebra already. That is, even if we were to chiralize an arbitrary Lie algebra $\mathfrak{g}$ and some arbitrary $\mathfrak{g}$-module $\mathfrak{v}$, assembled into the semidirect sum Lie algebra $\mathfrak{g} \oplus \mathfrak{v}$, the resulting vertex algebra would \textit{almost} contain the symmetric or antisymmetric algebras generated by $\mathfrak{v}$, and the coefficient ring of $\mathfrak{g}$ would \textit{almost} include $S\mathfrak{v} \otimes \Lambda \mathfrak{v}$.

Very little effort is needed to eliminate the recurring word ``almost" in the preceding paragraph. There is an obvious set of relations that must be added to the vertex algebra so that it possesses the souped-up Lie algebra structure. Including these relations is a crucial choice, not without consequences, and marks the departure of our construction from the construction of the chiral de Rham complex. As a result, our construction becomes a replacement for the chiral de Rham complex, rather than an extension of it.

The constituent vertex algebras in the chiral de Rham sheaf do not have the structure of a souped-up Lie algebra, which this author deems unsatisfactory. While the chiral de Rham sheaf appropriately transfers the action of a vector field $X$ on a smooth function $f$ to the vertex algebra (specifically, $Xf = X \circ_0f$), it does not transfer the \textit{product} of $f$ with a vector field $X$ to the vertex algebra (so $fX \neq f\circ_{-1}X$ in general). Our amendment is to impose the relation $fX = f\circ_{-1}X$. However, this does not come without a trade-off. To wit, although we are able to prove in some specific cases that a souped-up Lie algebra injects into the vertex algebra it generates, this injectivity is only \textit{conjectured} to hold for arbitrary souped-up Lie algebras, including for $\mathfrak{s} (\mathfrak{a} \otimes \mathfrak{s} \mathfrak{X}_{\nabla}, \mathfrak{a})$.

\subsubsection{Relation to other vertex algebra bundles}

Ours is not the first discussion of vertex algebra bundles. There are two others of which the author is aware, and both of those are substantially different from the one at hand. The one more closely related is that of \cite{DLMZ-2004}, who have defined a notion of a $K$-theory associated to vertex algebras. Their idea is to begin with the concept of a $V$-bundle on $M$, which is akin to a $G$-bundle but with the group $G$ replaced by a vertex algebra $V$. The difference between their bundle and ours is that in ours, not only is the structure group $G$ (whose action on the fibers is \textit{vertical}) replaced by a vertex algebra, but the fibers themselves are represented in this vertex algebra. And still furthermore, the Lie algebra of vector fields on $M$ (whose action on the fibers is \textit{horizontal}) is part of our vertex algebra. This final aspect (chiralizing the Lie algebra of vector fields) is what likens our construction to the chiral de Rham complex, and makes the situation both more intricate and robust than the $V$-bundles of \cite{DLMZ-2004}.

There is a second notion of a vertex algebra bundle that predates even the one already discussed. Frenkel and Ben-Zvi \cite{FBZ-2004} speak of a vertex algebra bundle over the base space $\C$ (some Riemann surface) in the quantum field theoretic point of view of conformal field theory, rather than over the target space $M$ as is done in our situation. Thus our bundles are placed in a very different context from theirs.

\subsection{Organization of this manuscript}

In Chapter \ref{sec:VA}, we give a definition of a vertex algebra, portrayed as the quotient of an algebra $F$ with an infinity of totally nonassociative noncommutative products $\circ_n$. The usual relations of a vertex algebra are then represented by an ideal $I \subset F$. This is accompanied by a change in philosophy regarding the nonnegative products $\circ_n$, $n\geq 0$. Whereas these products are typically chosen in accordance with the Ward identity (which relates these products to the transformation properties of the underlying classical geometry; see \cite{Polchinski-1998}), we now shed this obligation and permit these particular products to be specified arbitrarily, or even not at all.

Perhaps this new stance distances vertex algebras from their origins in string theory, but the result is a new object of interest in mathematics.

In Chapter \ref{sec:soupedupLA}, we define a new classical object, the \textit{souped-up Lie algebra}, which combines a Lie algebra $\mathfrak{g}$ and a module $\mathfrak{a}$ (an associative commutative unital algebra) into the semidirect sum Lie algebra $(\mathfrak{a} \otimes \mathfrak{g}) \oplus \mathfrak{a}$. We have enlarged the ring of coefficients of $\mathfrak{g}$ to include the algebra $\mathfrak{a}$. The most important souped-up Lie algebra we construct is the one associated to a $G$vector bundle $E$ with connection $\nabla$, discussed above in \ref{sec:CDRrelationship}.

In Chapter \ref{sec:soupedupVA}, we construct a vertex algebra from a souped-up Lie algebra. For some basic souped-up Lie algebras $\mathfrak{s}$, we prove that the resulting vertex algebra contains $\mathfrak{s}$ as a subspace. However, the viability of this construction for a \textit{generic} souped-up Lie algebra is ultimately a conjecture; the resulting vertex algebra might be zero (and thus not a vertex algebra at all). That said, we are confident that the resulting vertex algebra contains the entire generating souped-up Lie algebra. We prove some necessary conditions for these conjectures to be true.

The importance of Chapter \ref{sec:uVecttoVA} is the definition of the functor $\mathcal{V}$ from souped-up Lie algebras to vertex algebras. This functor will be used to associate a sheaf of souped-up Lie algebras to a sheaf of vertex algebras.

Next, in Chapter \ref{sec:sheaf}, we construct a sheaf of vertex algebras from a sheaf of souped-up Lie algebras.

Finally, in Chapter \ref{sec:chiralDG}, we apply this construction to the sheaf of souped-up Lie algebras associated to $(E,\nabla)$, giving us a sheaf of vertex algebras containing all the geometric and algebraic information from $(E,\nabla)$.

\subsubsection{A note on superspaces and supercommutativity}
Throughout this manuscript, we will assume that every vector space is in fact a \textit{superspace}, meaning it has a ${\Z}_2$-grading. Of course, it is certainly possible that the odd component is trivial. Furthermore, we will use \textit{commutative} to mean \textit{supercommutative}. Thus a commutative ring has a ${\Z}_2$-grading in which odd elements $x,y$ satisfy $xy + yx = 0$. Similarly, in an abelian Lie algebra $\mathfrak{g}$, the associative multiplication on the universal enveloping algebra for $\mathfrak{s}$ is given by $[x,y] = xy + yx$ for odd elements $x,y$.

Despite this convention, for the sake of clarity, we will still occasionally apply the prefix \textit{super-} when labeling certain vector spaces or products.

We will make explicit reference to the parity $p$ of an element only in Definition \ref{def:VA} of a vertex algebra ideal. After that, all definitions and results hold for both even and odd elements, but to simplify the presentation we will only prove them for the even case.

\section{Vertex algebras}
\label{sec:VA}
The notion of a vertex algebra was pinned down by Borcherds in 1986 \cite{borcherds-1986} to abstract the quantum fields prevalent in string theory. Although the axioms for a vertex algebra have been expressed in several different ways since their advent (see \cite{borcherds-1986}, \cite{kac-1997}), we have chosen to recast the definition yet again with a mind to the work that appears in later chapters. Because ours is not the orthodox definition, we will prove that it is equivalent to the original definition in Proposition \ref{prop:B-equivalence}.

\subsection{Definition of a vertex algebra as a quotient of an infinite free algebra}
We begin with an important preliminary definition of an \textit{infinite free algebra}. In our rendition, a vertex algebra is some quotient of this algebra.

We will often use the following notion of a \textit{unital vector space}.
\begin{defn}[Unital vector space]
A \textit{unital vector space} is a vector space with a distinguished vector $\1$.
\end{defn}

\begin{defn}[Infinite free algebra $F$]
\label{def:IFA}
An \textit{infinite free algebra} is a unital vector space that is closed under an infinite number of totally nonassociative noncommutative linear products $\circ_n$, $n\in {\Z}$. We will generally abbreviate the product $x \circ_n y$ by $x_ny$.

An \textit{ideal} in $F$ is closed with respect to \textit{all} products $\circ_n$.
\end{defn}

The element $\1$ is \textit{not} an identity for the products $\circ_n$ in general. However, upon taking a quotient of $F$ in the upcoming definition of a vertex algebra, $\1$ will become the two-sided identity for the product $\circ_{-1}$.

Every element in an infinite free algebra $F$ can be represented as a sum of full binary rooted trees in which each leaf represents another element of $F$, and each of the remaining nodes is labeled by an integer. As we have defined $F$, it is possible that a tree may be infinitely deep. This awkwardness is removed by assuming $F$ is generated by some unital vector space $\mathfrak{v}$, which we put in a definition.

\begin{defn}[Infinite free algebra $F(\mathfrak{v})$]
\label{def:generatedIFA}
We write $F(\mathfrak{v})$ when an infinite free algebra $F$ is generated by a unital vector space $\mathfrak{v}$. That is, $F(\mathfrak{v})$ is the closure of $\mathfrak{v}$ under the products $\circ_n$ and addition.
\end{defn}

Any additional structure on $\mathfrak{v}$, beyond its unital vector space structure, is forgotten in $F(\mathfrak{v})$.

In a \textit{generated} infinite free algebra $F(\mathfrak{v})$, one can now identify when a leaf node of some element $x$ is terminal: the leaf node is an element of $\mathfrak{v}$. We thus have the notion of a \textit{monomial} in $F(\mathfrak{v})$, which is an element that can be written as a product of elements of $\mathfrak{v}$, or can be depicted as a single tree in which every leaf node is in $\mathfrak{v}$.

While a vertex algebra may defined as a particular quotient of an infinite free algebra $F$ coming from Definition \ref{def:IFA}, virtually all examples in the literature, including those of importance in the current paper, are quotients of \textit{generated} infinite free algebras $F(\mathfrak{v})$.

\begin{note}[Multiples, factors, $\sembrack{\cdot}$, and $D$]
\label{note:notation}
In $F$ and $F(\mathfrak{v})$, we can speak of \textit{multiples} and \textit{factors}. When we refer to a \textit{multiple} $y$ of some element $x$, we mean that $y$ is obtained from $x$ by a sequence products on the left and right. In this case we also say that $x$ is a \textit{factor} of $y$. Equivalently, $x$ is a subtree of $y$.

We will often use the notation $y\sembrack{x}$ to denote the dependence of an element $y$ on another element $x$. In particular, $y\sembrack{x}$ means $x$ is a factor of some term in $y$.

We will follow the convention that the particular product $\circ_{-2}\1$ is denoted $D$, so $Dx \equiv x_{-2}\1$. We will see shortly that $D$ is a derivation over all products in a vertex algebra.
\end{note}

We now define a vertex algebra as a particular quotient of an infinite free algebra $F(\mathfrak{v})$.

\begin{defn}[Vertex algebra ideal $I(\mathfrak{v})$ and vertex algebra $V(\mathfrak{v})$]
\label{def:VA}
Let $F(\mathfrak{v})$ be an infinite free algebra, and fix some function $N(u,v) \geq 0$ on $\mathfrak{v} \times \mathfrak{v}$. The \textit{vertex algebra ideal} $I(\mathfrak{v})$ is the ideal generated by the following sets:
\begin{description}
\item[identity:] $\textbf{i}\sembrack{x;n} \triangleq \1_nx - \delta_{n+1}x$
\item[locality:] \hfill \\
$\textbf{c}\sembrack{u,v;n} \triangleq u_nv \quad \text{for all } n \geq N(u,v)$
\item[derivation:] \hfill \\
$\textbf{d}\sembrack{x,y;n} \triangleq D(x_ny) - (Dx)_ny - x_n(Dy)$ \\
$\textbf{e}\sembrack{x,y;n} \triangleq (Dx)_ny + nx_{n-1}y \quad \text{for all } n \in {\Z}$
\item[quasi-commutativity:] \hfill \\
$\textbf{qc}\sembrack{x,y;n} \triangleq x_ny + (-1)^{p(x) p(y)}\sum_{k\geq 0} \frac{(-1)^{n+k}}{k!} D^k(y_{n+k}x)$
\item[quasi-associativity:] \hfill \\
$\textbf{qa}\sembrack{x,y,z;m,n} \\ \triangleq (x_my)_nz - \sum_{k\geq 0}\binom{m}{k}(-1)^k \left(x_{m-k}(y_{n+k}z) - (-1)^{m+ p(x) p(y)} y_{m+n-k}(x_kz) \right)$
\end{description}
for $u,v\in \mathfrak{v}$ and $x,y,z \in F(\mathfrak{v})$, for $m,n=-1$, and for all $m,n \gg 0$ unless otherwise indicated. The binomial coefficient $\binom{m}{k}$ is defined as usual for all $m\in {\Z}$ and $k \in {\Z}_{\geq 0}$, and extends to all $k\in {\Z}$ by setting $\binom{m}{k} = 0$ if $k<0$. $p$ is the parity induced by the ${\Z}_2$-grading of $\mathfrak{v}$.

A vertex algebra $V(\mathfrak{v})$ is the quotient $F(\mathfrak{v})/J$, where $J$ any proper ideal containing the vertex algebra ideal $I(\mathfrak{v})$.
\end{defn}

By linearity, we may assume that the variables $x,y,z$ represent \textit{monomials}. This assumption simplifies future arguments.

Although we will not be needing it much, we give the definition of a vertex algebra in the more general case that $F$ is not assumed to be generated by $\mathfrak{v}$.

\begin{defn}[Vertex algebra ideal $I$ and vertex algebra $V$]
\label{def:generalVA}
In the case that $F$ is not generated by some $\mathfrak{v}$, we define the vertex algebra ideal $I$ more generally by replacing the function $N(u,v)$ with a function $N(x,y) \geq 0$ for all $x,y \in F$. Then a vertex algebra is the quotient $V = F/J$, where $J$ is any proper ideal containing $I$.
\end{defn}

Upon forming the quotient $V$ (or $V(\mathfrak{v})$), the elements $\textbf{i}\sembrack{x;n}$ imply that $\1$ is a \textit{left} identity for the product $\circ_{-1}$ on $V$. Making use of the quasi-commutativity relation $\textbf{qc}\sembrack{x,\1;-1}$ on $V$, we see that $\1$ is a \textit{right} identity for $\circ_{-1}$ as well. The elements $\textbf{c}\sembrack{x,y;n}$ ensure that the tail ends of the series appearing in $\textbf{qc}\sembrack{x,y;n}$ and $\textbf{qa}\sembrack{x,y,z;m,n}$ get killed in $V$, so that the series are in fact convergent. The elements $\textbf{d}\sembrack{x,y;n}$ express that $D \equiv \circ_{-2}\1$ is a derivation of all products in $V$. The elements $\textbf{qc}\sembrack{x,y;n}$ and $\textbf{qa}\sembrack{x,y,z;m,n}$ express the extent to which each product is not commutative or associative on $V$, as seen by comparing the leading term with the first term of the summation.

\begin{note}
We will write $F$ and $V$ when we are making no assumptions about the existence of an underlying vector space $\mathfrak{v}$.

We will use the letters $u,v$ to denote elements of $\mathfrak{v}$ and their inclusions in $F(\mathfrak{v})$; the letters $x,y$ denote either general elements of some infinite free algebra.  Thus when we refer to $N(u,v)$, we mean a function $N$ defined on $\mathfrak{v} \times \mathfrak{v}$, whereas $N(x,y)$ denotes a function defined on all of $F\times F$ or $F(\mathfrak{v}) \times F(\mathfrak{v})$. Also, an element $x\in F, F(\mathfrak{v})$ will also be used to denote its equivalence class $x+J \in V,V(\mathfrak{v})$.

From now on, to ease the notation, we will assume all elements are even, so that the parity $p(x)=0$ for all $x\in F, F(\mathfrak{v})$. Even so, with an appropriate adjustment of signs, every subsequent definition and result still applies to both even and odd elements.
\end{note}

\begin{rmk}
The impact of the $u,v$-dependence of the integer $N(u,v)$ in the second set of generators $\textbf{c}\sembrack{u,v;n}$ is subtle but important. We are \textit{not} guaranteed that there is a single number $n$ such that $u_nv \in I(\mathfrak{v})$ for all $u,v \in \mathfrak{v}$, but rather that for each pair $u,v$ there is a number $N(u,v)$ such that $u_nv \in I(\mathfrak{v})$ for all $n \geq N(u,v)$. In contrast, in all of the other sets of generators, the choice of $n$ is independent of $x$ and $y$. It is the content of Dong's Lemma (proved in Proposition \ref{prop:dong}) that in fact for any pair of elements $x,y \in F(\mathfrak{v})$ there is some number $N(x,y)\geq 0$ such that $x_ny \in I(\mathfrak{v})$ for all $n\geq N(x,y)$.

The function $N(u,v)$ has not been specified, meaning that the notation ideal $I(\mathfrak{v})$ has some ambiguity. Thus when we speak of \textit{the} vertex algebra ideal $I(\mathfrak{v})$, we really mean \textit{some} vertex algebra ideal $I(\mathfrak{v})$ for which $N(u,v)$ has been specified.

In contrast, the lack of specification of the integers $m$ and $n$ is harmless. As we will see in Proposition \ref{prop:B-equivalence}, the inclusion of any generator for $n=-1$ and $n\gg 0$ implies the inclusion of that family of generators for all $n \in {\Z}$.
\end{rmk}

\subsection{Equivalence to Borcherds' definition}

For this section, we increase our scope to $F$ and $V$ as in Definition \ref{def:generalVA}, not assumed to be generated by $\mathfrak{v}$. Thus $I$ includes the larger set of generators $\textbf{c}\sembrack{x,y;n} = x_ny$ for $n \geq N(x,y)$ for some function $N(x,y) \geq 0$ on $F\times F$.

The key difference between our definition and others is the range of the integers $m, n$ in the products $\circ_m, \circ_n$ appearing in the generators of $I$. In our definition, the number of generators is severely reduced in that generally we have only $m,n=-1$ and all $m,n \gg 0$, whereas in the standard definitions $m$ and $n$ usually range over all of ${\Z}$. Our only set of generators in which the product $\circ_n$ is indexed by all $n \in {\Z}$ is $\textbf{e}\sembrack{x,y;n}$. As we will see in the upcoming Proposition \ref{prop:B-equivalence}, with the aid of this particular set of generators, all standard generators can be recovered.

\begin{prop}
\label{prop:B-equivalence}
Definition \ref{def:generalVA} of a vertex algebra is equivalent to Borcherds' original definition.
\end{prop}

\begin{proof}
Borcherds' relations in a vertex algebra are almost the same as our generators for $I$. On one hand, his versions of $\textbf{i}\sembrack{x;n}$, $\textbf{d}\sembrack{x,y;n}$, $\textbf{qc}\sembrack{x,y;n}$, and $\textbf{qa}\sembrack{x,y,z;m,n}$ hold for \textit{all} $m,n \in {\Z}$. But offsetting these extra relations, his set lacks our $\textbf{e}\sembrack{x,y;n}$.

To see that Borcherds' relations imply ours, one can check explicitly that $\textbf{e}\sembrack{x,y;n}$ is the sum of elements from Borcherds' relations:
\begin{align*}
\textbf{e}\sembrack{x,y;n} &= \textbf{qa}\sembrack{x,\1,y;-2,n} \\
&\quad + \sum_{k\geq 0}\binom{-2}{k}(-1)^k \left(x_{-2-k}\textbf{i}\sembrack{y;n+k} - \textbf{i}\sembrack{x_ky;-2+n-k} \right).
\end{align*}

For the other direction, we must show that $I$ contains the elements $\textbf{d}\sembrack{x,y;n}$, $\textbf{qc}\sembrack{x,y;n}$, and $\textbf{qa}\sembrack{x,y,z;m,n}$ for all $m,n \in {\Z}$. The argument is a downward induction.

We begin by showing that $\textbf{d}\sembrack{x,y;n-1} \in I$ whenever $\textbf{d}\sembrack{x,y;n} \in I$. As a base case, we are given that $\textbf{d}\sembrack{x,y;n} \in I$ for all $n \gg 0$. We have the readily checked identity
\begin{align*}
n\textbf{d}\sembrack{x,y;n-1} = D\textbf{e}\sembrack{x,y;n} - \textbf{d}\sembrack{Dx,y;n} - \textbf{e}\sembrack{Dx,y;n} - \textbf{e}\sembrack{x,Dy;n}.
\end{align*}
We see that on the left the index of $\textbf{d}$ is $n-1$ whereas on the right the index of $\textbf{d}$ is $n$, so $I$ contains $\textbf{d}\sembrack{x,y;n-1}$ whenever it contains $\textbf{d}\sembrack{x,y;n}$. This is precisely the (downward) inductive argument. This breaks down when $n=0$ since the left side vanishes. The base case is re-founded with the inclusion of $\textbf{d}\sembrack{x,y;-1}$ in our set of generators. In conclusion, $\textbf{d}\sembrack{x,y;n} \in I$ for all $n \in {\Z}$.

For the rest of the proof we abbreviate $\textbf{d}\sembrack{x,y;n} + \textbf{e}\sembrack{x,y;n}$ as $\textbf{f}\sembrack{x,y;n}$ for all $n \in {\Z}$.

To show that $\textbf{i}\sembrack{x;n}, \textbf{qc}\sembrack{xy;n}, \textbf{qa}\sembrack{xy,z;m,n} \in I$ for all $n \in \Z$, $I$ already includes the base cases for all $m,n \gg 0$ and $m,n=-1$. The separate inductive arguments are then given by the equalities
\begin{align*}
n\textbf{i}\sembrack{x;n-1} &= \textbf{i}\sembrack{Dx;n} - D\textbf{i}\sembrack{\1,x;n} + \textbf{f}\sembrack{\1,x;n}, \\
n\textbf{qc}\sembrack{x,y;n-1} &= -\textbf{qc}\sembrack{Dx,y;n} + \textbf{e}\sembrack{x,y;n} - \sum_{k\geq 0} \frac{(-1)^{n+k}}{k!} D^k\textbf{f}\sembrack{y,x;n+k},
\end{align*}

\begin{align*}
m \textbf{qa}\sembrack{x,y,z;&m-1,n} = -\textbf{qa}\sembrack{Dx,y,z;m,n} + \textbf{e}\sembrack{x,y;m}_nz \\
& - \sum_{k\geq 0}\binom{m}{k}(-1)^k \left(\textbf{e}\sembrack{x,y_{n+k}z;m-k} - (-1)^m y_{m+n-k} \textbf{e}\sembrack{x,z;k} \right),
\end{align*}
and
\begin{align*}
n \textbf{qa}\sembrack{x,y,z;&m,n-1} = -D\textbf{qa}\sembrack{x,y,z;m,n} + \textbf{qa}\sembrack{x,y,Dz;m,n} + \textbf{f}\sembrack{x_my,z;n} \\
& - \sum_{k \geq 0}\binom{m}{k}(-1)^k \left( \textbf{f}\sembrack{x,y_{n+k}z;m-k} - (-1)^m \textbf{f}\sembrack{y,x_kz;m+n-k} \right) \\
& - \sum_{k \geq 0}\binom{m}{k}(-1)^k \left( x_{m-k} \textbf{f}\sembrack{y,z;n+k} - (-1)^m y_{m+n-k} \textbf{f}\sembrack{x,z;k}\right).
\end{align*}
\end{proof}

\begin{rmk}
\label{rmk:fields}
Other axioms for a vertex algebra feature a collection of \textit{quantum fields} $x(\zeta) \in (\End V)[[\zeta, \zeta^{-1}]]$ (formal power series in $\zeta$ and $\zeta^{-1}$ with coefficients in $\End V$), and an even endomorphism $D$ of $V$ satisfying $[D,x(\zeta)] = \frac{d}{d\zeta} x(\zeta)$. Our formulation captures this data via the definitions
\begin{align*}
x(\zeta) &\triangleq \sum_{n\in \Z}\frac{x \circ_n}{\zeta^{n+1}} \\
D &\triangleq \circ_{-2}\1.
\end{align*}
\end{rmk}

\subsection{Decompositions of \texorpdfstring{$F(\mathfrak{v})$}{F(v)}}
\label{sec:decompositions}

$F(\mathfrak{v})$ has some useful linear decompositions. One of the decompositions we discuss, the \textit{degree grading}, descends to a decomposition of the quotient $V(\mathfrak{v}) = F(\mathfrak{v}) / I(\mathfrak{v})$, precisely because each generator of $I(\mathfrak{v})$ has homogenous degree.

First we will discuss a set of decompositions of $F(\mathfrak{v})$ that are ordered by coarseness. To begin, we define a \textit{monomial} as any element that can be written as a single (full binary rooted) tree whose leaves are elements of $\mathfrak{v}$. Then a very fine decomposition is of $F(\mathfrak{v})$ into a direct sum of subspaces spanned by single monomials. As examples, $\{v,\, v_2w,\, (u_{-2}v)_1w\}$ are monomials provided $u,v,w \in \mathfrak{v}$. The element $u_1v + w_1v$ is also a monomial because it can be written as $(u+w)_1v$. Two monomials are in the same summand precisely when corresponding leaf nodes are scalar multiples of one another.

Slightly coarser, we can group monomials by \textit{product shape}, which pays attention no not only to be the shape of the tree, but also to the integer at each node. The product shape captures the sequence of products $\circ_n$ used to construct this element. Examples of elements of homogeneous product shape are $\{v_2w,\, u_5(v_2w) + x_5(y_2u),\, (u+v)_{-1}w + x_{-1}y\}$, where the variables represent monomials. Note that any monomial has homogeneous product shape. An element of homogeneous product shape can be represented by replacing all leaves with asterisks, in which case our three examples are now written $\{\ast_2\ast,\, \ast_5(\ast_2\ast) + \ast_5(\ast_2\ast),\, \ast_{-1}\ast + \ast_{-1}\ast\}$.

Coarser still is the \textit{shape} of the tree alone, without regard to the integers at each node. In this case, we can ignore the products $\circ_n$ altogether, and simply group leaves with parentheses. Then the previous example becomes $\{\ast\ast,\, \ast(\ast\ast) + \ast (\ast\ast),\, \ast\ast + \ast\ast\}$.

And coarsest of all, we have the \textit{length} of a monomial, which is simply its number of leaves $\ast$. Altogether, we have the ordering, from coarse to fine,
\[
\begin{tikzpicture}
\matrix(m)[matrix of math nodes,
row sep=2em, column sep=3em,
text height=1.5ex, text depth=0.25ex]
{ \text{length} &  \text{shape} &  \text{product shape} &  \text{monomial} \\};
\path[->]
	(m-1-1) edge (m-1-2)
	(m-1-2) edge (m-1-3)
	(m-1-3) edge (m-1-4);
\end{tikzpicture}
\]

Note that the subspace $\mathfrak{v} \subset F(\mathfrak{v})$ is precisely the subspace of length 1. Thus in the decomposition by length, the subspaces of lengths 2 or greater are orthogonal to $\mathfrak{v}$.

We now present a degree grading that leads to a decomposition descending nicely onto the quotient $V(\mathfrak{v})$. Equivalently, the vertex algebra ideal $I(\mathfrak{v})$ decomposes according to the degree grading.

Given a direct sum decomposition of $\mathfrak{v}$, and assigning each summand a weight (not necessarily in ${\Z}$), there is a number of ways a degree grading $|\cdot|$ is induced on $F(\mathfrak{v})$. Insisting that the element $\1 \in \mathfrak{v}$ has degree 0, and insisting that the grading be additive over the products $\circ_n$, meaning
\[ |x_ny| = |x| + |\circ_n| + |y|, \]
then such a grading will descend nicely to the vertex algebra $F(\mathfrak{v})/I(\mathfrak{v})$ only when the degree of the product $\circ_n$ is given by
\[ |\circ_n| \triangleq - n - 1. \]
Indeed, in this case it is easily verified that every generator of $I$ has homogeneous degree.

The decomposition of $F(\mathfrak{v})$ by degree is neither finer nor coarser than the decompositions above, but a common refinement can be found, giving the partial ordering diagram
\[
\begin{tikzpicture}
\matrix(m)[matrix of math nodes,
row sep=2em, column sep=1.3em,
text height=1.5ex, text depth=0.25ex]
{ & \text{length} &  \text{shape} &  \text{prod. shape} &  \text{monomial} \\ {|\cdot|} & \text{length} \cap |\cdot| &  \text{shape} \cap |\cdot| &  \text{prod. shape} \cap |\cdot| &  \text{monomial} \cap |\cdot| \\};
\path[->]
	(m-1-2) edge (m-1-3)
	(m-1-3) edge (m-1-4)
	(m-1-4) edge (m-1-5)
	(m-2-1) edge (m-2-2)
	(m-2-2) edge (m-2-3)
	(m-2-3) edge (m-2-4)
	(m-2-4) edge (m-2-5)
	(m-1-2) edge (m-2-2)
	(m-1-3) edge (m-2-3)
	(m-1-4) edge (m-2-4)
	(m-1-5) edge (m-2-5);
\end{tikzpicture}
\]

The generators of the ideal $I(\mathfrak{v})$ are homogeneous in none of these decompositions except for the degree grading $|\cdot|$. As a result, only the degree grading descends onto the vertex algebra $V(\mathfrak{v})$. Despite this fact, the we are still able to speak of \textit{monomials} and \textit{lengths} in $V(\mathfrak{v})$.

\begin{defn}[Monomials and length in $V(\mathfrak{v})$]
\label{def:VAmonomiallength}
Let $V(\mathfrak{v}) = F(\mathfrak{v}) / J$ for some ideal $J \supset I(\mathfrak{v})$. We say an element $x \in V(\mathfrak{v})$ is a \textit{monomial} if the coset $x + J$ contains a monomial in $F(\mathfrak{v})$. The \textit{length} of $x$ is the length, in $F(\mathfrak{v})$, of the shortest element in the coset $x + J$.
\end{defn}

\subsection{Injectivity of \texorpdfstring{$\mathfrak{v}$ into $V(\mathfrak{v}) = F(\mathfrak{v})/I(\mathfrak{v})$}{\textit{v} into \textit{V(v) = F(v)/I(v)}}}

A reassuring feature of a vertex algebra generated by a unital vector space $\mathfrak{v}$ is that $\mathfrak{v}$, considered as a subspace of $F(\mathfrak{v})$, survives intact upon taking the quotient by $I(\mathfrak{v})$. We will delay the proof of this statement until Chapter \ref{sec:soupedupVA}, where it will be a corollary to the slightly stronger Theorem \ref{thm:Ainjectivity}. In that theorem, we will prove that $\mathfrak{v}$ survives for a particular enlargement of the ideal $I(\mathfrak{v})$.

For now, we state this fact as a theorem without proof.
\begin{thm}
\label{thm:Vinjectivity}
The map taking $\mathfrak{v}$ to its image in $V(\mathfrak{v}) = F(\mathfrak{v})/I(\mathfrak{v})$ is a monomorphism. Equivalently, $I(\mathfrak{v}) \cap \mathfrak{v} = \{0\}$.
\end{thm}

This injectivity is independent of choice of the function $N(u,v)$. In particular, we may choose the most constrictive function possible, setting $N(u,v) \equiv 0$, which has the effect of enlarging $I(\mathfrak{v})$ and shrinking $V(\mathfrak{v})$ simultaneously. In particular, this choice implies that in $V(\mathfrak{v})$, $u_nv = 0$ for all $n\geq 0$ and all $u,v \in \mathfrak{v}$.

To keep our perspective, we point out that this injectivity could be destroyed if we enlarge $I(\mathfrak{v})$ further in other ways. We could easily collapse a part of $\mathfrak{v}$ by replacing $I(\mathfrak{v})$ with a larger ideal $J \supset I(\mathfrak{v})$ that, say, includes some element $v \in \mathfrak{v}$.

\subsection{Commutator formula and Dong's Lemma}
In this section, for completeness, we present the long-established commutator formula between the \textit{mode operators} $x_m$ and $y_n$ on general $F$ and $V$ (as in Definition \ref{def:generalVA}). This formula helps us prove Dong's Lemma \ref{prop:dong} (originally proven in \cite{Li:1994sp}; see \cite{kac-1997} for another version), which is essential for the convergence of the quasi-commutativity and quasi-associativity relations in a vertex algebra $V(\mathfrak{v})$. In particular, it guarantees that vanishing of the products $u_nv$ for $u,v \in \mathfrak{v}$ and sufficiently large $n$ implies the vanishing of \textit{all} products $x_ny$ for all $x,y \in V(\mathfrak{v})$ for sufficiently large $n$.

\begin{prop}
\label{prop:commutator}
In a vertex algebra $V$, we have the commutator formula
\[ [x_m,y_n]z = \sum_{k\geq 0} \binom{m}{k} (x_ky)_{m+n-k}z.\]
\end{prop}

\begin{proof}
It is straightforward to verify the identity
\begin{align*}
[x_m,y_n]z - \sum_{k\geq 0} \binom{m}{k}& (x_ky)_{m+n-k}z = \textbf{qa}\sembrack{x,y,z;-1,-1} - \textbf{qa}\sembrack{y,x,z;-1,-1} \\
& - \textbf{qc}\sembrack{y,x;-1} + \sum_{0\leq i \leq j} \frac{(-1)^{j+1} i!}{(j+1)!} \textbf{e}\sembrack{D^{j-i}(x_jy),z; -1-i}.
\end{align*}
The proposition then follows since all terms on the right side are in $I$, and therefore so is the left side, which is our desired relation.
\end{proof}

\begin{lem}
\label{lem:Nsymmetric}
$N(x,y)$ is symmetric. That is, if $x_ny \in I$ for all $n\geq N(x,y)$, then also $y_nx \in I$ for $n\geq N(x,y)$.
\end{lem}

\begin{proof}
Assuming that $x_ny \in I$ for $n \geq N(x,y)$, we have
\[ y_nx = \textbf{qc}\sembrack{y,x;n} - \sum_{k\geq 0} \frac{(-1)^{n+k}}{k!} D^k(x_{n+k}y). \]
Every term on the right is in $I$, thus so is $y_nx$.
\end{proof}

We now prove Dong's Lemma for a vertex algebra $V(\mathfrak{v})$ generated by a unital vector space $\mathfrak{v}$.

\begin{prop}[Dong's Lemma]
\label{prop:dong}
For any pair of elements $x,y \in F(\mathfrak{v})$, there exists some $N(x,y)\geq 0$ such that $x_ny \in I(\mathfrak{v})$ for all $n\geq N(x,y)$.
\end{prop}

\begin{proof}
By linearity, we may assume that $x$ and $y$ are monomials. We will perform an induction on the length of $x_ny$.

The base case, length 2, is given to us by the inclusion $\textbf{c}\sembrack{u,v;n} \in I(\mathfrak{v})$ for all $n \geq N(u,v)$ and $u,v \in \mathfrak{v}$.

Now assume that the lemma is true for all monomials of length $p$. We wish to show that for a monomial $w_nz$ having length $p+1$ there exists some $N$ such that $w_nz \in I(\mathfrak{v})$ for all $n\geq N$. At least one of $w$ and $z$ necessarily factors further. By Lemma \ref{lem:Nsymmetric}, we may assume that $w$ factors, so that $w_nz = (x_ry)_nz$. By hypothesis, since any pair among $x,y$, and $z$ has length not exceeding $p$, there exists a sufficiently large $M$ such that all products $x_my, y_mz, x_mz$, for $m \geq M$, are in $I(\mathfrak{v})$. ($M$ may be the maximum of $N(x,y), N(x,z)$, and $N(y,z)$.)

We will proceed in two steps: first we will prove the proposition for $r\geq 0$, and then for $r<0$.

Consider the commutator formula from Proposition \ref{prop:commutator} in the case $m \geq 2M, n=m-M$. The two terms $x_m(y_{m-M}z)$ and $y_{m-M}(x_mz)$ on the left side both vanish in $V(\mathfrak{v})$ since (as elements of $F(\mathfrak{v})$) they are in $I(\mathfrak{v})$ by hypothesis. Similarly, on the right side, all terms with $k\geq M$ vanish, leaving altogether
\[ \sum_{k=0}^{M-1} \binom{m}{k} (x_ky)_{2m-M-k}z = 0.\]
Consider as well the cases $(m-1,m-M+1), (m-2, m-M+2), \ldots, (m-M, m)$, giving the equations
\[ \sum_{k=0}^{M-1} \binom{m-1}{k} (x_ky)_{2m-M-k}z = 0,\,\ldots,\, \sum_{k=0}^{M-1} \binom{m-M}{k} (x_ky)_{2m-M-k}z = 0.\]
Altogether we have a system of $M+1$ homogeneous equations in $M$ unknowns, and the binomial coefficients force the solution to be trivial:
\[ (x_0y)_{2m-M}z = (x_1y)_{2m-M-1}z = \cdots = (x_{M-1}y)_{2m-2M+1}z = 0.\]
Equivalently, each monomial above is in the ideal $I(\mathfrak{v})$. This holds for all $m \geq 2M$, so for a monomial of length $p+1$ of the form $(x_ry)_nz$ with $r\geq 0$, the proposition is satisfied by the choice $n \geq N = 3M - r$. By Lemma \ref{lem:Nsymmetric}, this also holds for elements $z_n(x_ry)$, $r\geq 0$.

To prove the lemma for length $p+1$ elements of the form $(x_ry)_nz$ for $r<0$, we note that when $n \geq N = 3M-r$, then every term on the right side of the equation
\[ (x_ry)_nz  = \textbf{qa}\sembrack{x,y,z;r,n} + \sum_{k\geq 0}\binom{r}{k}(-1)^k \left(x_{r-k}(y_{n+k}z) - (-1)^r y_{r+n-k}(x_kz) \right) \]
is in the ideal $I(\mathfrak{v})$. Indeed, the terms $x_{r-k}(y_{n+k}z)$ in the summation are in the ideal by hypothesis, since $n+k \geq M$; for the same reason, the terms $y_{r+n-k}(x_kz)$ for $k \geq M$ are in the ideal; and the terms $y_{r+n-k}(x_kz)$ for $0 \leq k < M$ are in the ideal in light of the result above for $r \geq 0$. Thus $(x_ry)_nz \in I(\mathfrak{v})$ for all $n \geq N = 3M-r$ when $r < 0$.

This concludes the proof of Dong's Lemma.
\end{proof}

\section{Souped-up Lie algebras}
\label{sec:soupedupLA}

In the well-known \textit{current algebra construction} (see \cite{kac-1997}, Chapter 2.5), a vertex algebra is generated from a Lie algebra $\mathfrak{g}$. This vertex algebra actually contains a one-dimensional central extension of $\mathfrak{g}$, with the Lie bracket played by the product $\circ_0$, and the central extension generated by the vacuum element $\1$.

With some additional simple relations (addressed next in Chapter \ref{sec:soupedupVA}), a vertex algebra may contain a much richer classical structure, combining a Lie algebra $\mathfrak{g}$ and a $\mathfrak{g}$-module $\mathfrak{a}$. In this chapter we define this structure, a \textit{souped-up Lie algebra}, and give several popular examples. Our final example will be a combination of the exterior algebra of differential forms $\Omega(M)$, the Lie algebra of contraction operators $\iota_X$ by vector fields $X \in \mathfrak{X}(M)$, sections of a vector bundle $E\to M$, the Lie algebra of infinitesimal gauge transformations on the space of sections, and a connection $\nabla$. This robust souped-up Lie algebra will be the basis of our final construction: a chiral vector bundle, in Chapter \ref{sec:chiralDG}.

\subsection{Definition of a souped-up Lie algebra}
We begin by defining a \textit{souped-up Lie algebra}.

\begin{defn}[Souped-up Lie algebra]
\label{def:soupedupLA}
Let $\mathfrak{a}$ be an associative commutative unital algebra and let $\mathfrak{g}$ be a Lie algebra such that
\begin{enumerate}
\item $\mathfrak{g}$ is a left $\mathfrak{a}$-module, this product written as juxtaposition.
\item $\mathfrak{a}$ is a two-sided $\mathfrak{g}$-module in which $\mathfrak{g}$ acts by derivations, written as the Lie bracket $[\cdot, \cdot]$.
\item The two products are compatible in that $[g,ah] = [g,a]h + a[g,h]$ for $a \in \mathfrak{a}$ and $g,h \in \mathfrak{g}$.
\end{enumerate}
Interpreting $\mathfrak{a}$ as an abelian Lie algebra, the \textit{souped-up Lie algebra} $\mathfrak{s}(\mathfrak{g}, \mathfrak{a})$ is defined as the semidirect sum Lie algebra $\mathfrak{g} \oplus \mathfrak{a}$ with the left multiplication by $\mathfrak{a}$.
\end{defn}

It is permissible that a souped-up Lie algebra has additional relations, and it will still be considered a souped-up Lie algebra as long as $\1\neq 0$.

\begin{rmk}
\label{rmk:soupedupLA}
It might be tempting to extend the left multiplication of $\mathfrak{a}$ on $\mathfrak{g}$ to a two-sided multiplication by symmetry, effectively adding the relation $s_{-1}a = as$, but as we will see below in Proposition \ref{prop:leftmult}, this relation is too strong to carry over into a vertex algebra, even though it poses no problems to the underlying souped-up Lie algebra.
\end{rmk}

\subsection{Examples of souped-up Lie algebras}
\label{sec:SLAexamples}

Souped-up Lie algebras include the following examples:

\begin{example}[Associative commutative unital algebra]
An associative commutative unital algebra $\mathfrak{a}$ is isomorphic to the souped-up Lie algebra $\mathfrak{s} (0, \mathfrak{a})$.
\end{example}

\begin{example}[Centrally extended Lie algebra]
Let $\mathfrak{g}$ be a Lie algebra, and let a new element $c$ generate a central extension of $\mathfrak{g}$. Letting $\mathfrak{a} = \spanof\{c\}$, we have a souped-up Lie algebra $\mathfrak{s}(\mathfrak{a} \otimes \mathfrak{g}, \mathfrak{a})$.
\end{example}

\begin{example}[Supersymmetric algebra of a Lie algebra module]
Let $\mathfrak{g}$ be a Lie algebra and $\mathfrak{v}$ a $\mathfrak{g}$-module. Let $\mathfrak{a}$ be the associative commutative unital algebra $S\mathfrak{v} \otimes \Lambda \mathfrak{v}$, which is also a module for $\mathfrak{g}$. We may convert $\mathfrak{g}$ to a left module for $\mathfrak{a}$ by extending its coefficients: $\mathfrak{a} \otimes \mathfrak{g}$. Altogether, we have the souped-up Lie algebra $\mathfrak{s}(\mathfrak{a} \otimes \mathfrak{g}, \mathfrak{a})$. The bracket on the Lie algebra underlying $\mathfrak{s}$ is given by
\[ [a\otimes s, b] \triangleq a \otimes [s,b]. \]
\end{example}

\begin{example}[Vector fields and smooth functions]
Let $\mathfrak{X}$ be the Lie algebra of vector fields on a smooth manifold $M$, and let $\Omega^0$ be the algebra of smooth functions on $M$. $\mathfrak{X}$ acts on $\Omega^0$ by the Lie derivative and $\mathfrak{X}$ is a left module for ${\Omega}^0$, the latter multiplication compatible with the Lie derivative action. Then $\mathfrak{s}(\mathfrak{X}, {\Omega}^0)$ is a souped-up Lie algebra.
\end{example}

\begin{example}[Lie superalgebra of vector fields and differential forms]
Let $\mathfrak{sX}$ be the Lie superalgebra of vector fields on $M$, and let $\Omega$ be the superalgebra of differential forms on $M$. $\mathfrak{sX}$ is the semidirect sum Lie algebra formed by $\mathfrak{X}$ acting on $\Omega$ as Lie derivatives ($(X,\omega) \mapsto \mathcal{L}_X \omega)$, and a second copy of $\mathfrak{X}$ acting as the \textit{odd} Lie algebra of contraction operators ($(Y,\omega) \mapsto \iota_Y \omega$). The Lie bracket between any two elements of  $\mathfrak{sX}$ is given by their commutator as operators on $\Omega$. Altogether, this leads to the bracket
\begin{align*}
[ \mathcal{L}_X, \mathcal{L}_Y] &= \mathcal{L}_{[X,Y]} \\
[\mathcal{L}_X, \iota_Y] &= \iota_{[X,Y]} \\
[\iota_X, \iota_Y] &= 0.
\end{align*}
We elevate this Lie algebra and its module to a souped-up Lie algebra by extending the coefficients of $\mathfrak{sX}$ to include $\Omega$, giving altogether $\mathfrak{s} ( \Omega \otimes_{\Omega^0} \mathfrak{s} \mathfrak{X}, \Omega)$.
\end{example}

\begin{example}[de Rham complex]
\label{ex:deRham}
We can add the exterior derivative $d$, regarded as an \textit{odd} element, to the Lie algebra $\mathfrak{s}\mathfrak{X}$ from the previous example, forming the Lie algebra $\mathfrak{s}\mathfrak{X}_d$. In fact, this Lie algebra is entirely generated by the elements $\{\iota_X \mid X\in \mathfrak{X}\}$ and $d$, since we have the odd bracket $[d, \iota_X] = \mathcal{L}_X$ (the famous Cartan formula for $\mathcal{L}_X$). Altogether $\mathfrak{s}\mathfrak{X}_d$ has the brackets defined on $\mathfrak{s}\mathfrak{X}$, plus
\begin{align*}
[\mathcal{L}_X, d] &= 0 \\
[d, \iota_X] &= \mathcal{L}_X \\
[d,d] &= 0.
\end{align*}
The associated souped-up Lie algebra is then  $\mathfrak{s} ( \Omega \otimes_{\Omega^0} \mathfrak{s} \mathfrak{X}_d, \Omega)$.
\end{example}

\begin{example}[Gauge Lie algebra and sections of a vector bundle]
\label{ex:gauge}
Consider a smooth vector bundle $E\to M$ with structure Lie group $G$. The group of gauge transformations $\mathcal{G}$ is a subbundle of $\Gamma (\Aut E)$. The Lie algebra $L\mathcal{G}$ to the gauge group is then a subbundle of $\Gamma (\End E)$. The action of $L\mathcal{G}$ on the space of sections $\Gamma E$ induces an action on the spaces of sections on the symmetric and antisymmetric tensor algebras $SE$ and $\Lambda E$. Then we can form the souped-up Lie algebra $\mathfrak{s}(\Gamma (SE \otimes \Lambda E) \otimes_{\Omega^0} L\mathcal{G}, \Gamma (SE \otimes \Lambda E))$.
\end{example}

\begin{example}[Connection Lie algebra and sections of a vector bundle]
\label{ex:connectionLA}
This final example is of primary importance to this paper. We combine Examples \ref{ex:deRham} and \ref{ex:gauge} to form what we call the \textit{connection Lie algebra} $\mathfrak{s} \mathfrak{X}_{\nabla}$.

Consider a $G$-vector bundle with connection $(E,\nabla)$. We designate as $\mathfrak{a}$ the algebra $\Omega \otimes \Gamma (SE \otimes \Lambda E)$, the tensor product of the exterior algebra of differential forms with the supersymmetric tensor algebra of $E$.

The Lie algebra $\mathfrak{g}$ is the subalgebra of $\End \mathfrak{a}$ generated by the contraction operators $\{\iota_X \mid X \in\mathfrak{X}\}$, the connection $\nabla$, and the gauge Lie algebra $L\mathcal{G}$ from the previous example. The first set of operators (all odd elements) acts on $\omega \otimes s \in \Omega \otimes \Gamma (SE \otimes \Lambda E)$ by performing as usual on the first factor and ignoring the second. The connection, also an odd element, acts as usual on the entire product, and an infinitesimal gauge transformation $A$, an even element, acts as usual on the second factor and ignores the first. Altogether, we have the actions (for $\omega$ with homogeneous degree $p$)
\begin{align*}
\iota_X (\omega \otimes s) &= \iota_X\omega \otimes s \\
\nabla (\omega \otimes s) &= d\omega \otimes s + (-1)^p \omega \otimes \nabla s\\
A (\omega \otimes s) &= \omega \otimes As.
\end{align*}
The Lie bracket between any pair is simply their commutator as operators on $\mathfrak{a}$. We denote this Lie algebra $\mathfrak{s}\mathfrak{X}_{\nabla}$.

It is straightforward to check that $[\iota_X,\iota_Y] = [\iota_X, A] = 0$. Of course, the bracket of $\nabla$ with itself is $[\nabla,\nabla] = 2\nabla^2$, twice the curvature 2-form.

Finally, we combine the algebra $\Omega \otimes \Gamma (SE \otimes \Lambda E)$ with the Lie algebra $\mathfrak{s} \mathfrak{X}_{\nabla}$ to form the souped-up Lie algebra
\[ \mathfrak{s} \equiv \mathfrak{s}( \Omega \otimes \Gamma (SE \otimes \Lambda E) \otimes \mathfrak{s} \mathfrak{X}_{\nabla}, \Omega \otimes \Gamma (SE \otimes \Lambda E)), \]
where all tensor products are over $\Omega^0$.

The souped-up connection Lie algebra contains all the standard algebraic and geometric information regarding the underlying manifold $M$. The subspace corresponding to $i=j=k=0$ in the first term in the direct sum decomposition
\begin{align}
\label{eq:soupconnLAdecomp}
&\mathfrak{s} = \bigoplus_{0 \leq i,j,k} \Omega^i\otimes \Gamma (S^jE) \otimes  \Gamma (\Lambda^k E) \otimes \mathfrak{s} \mathfrak{X}_{\nabla} \oplus \bigoplus_{0 \leq i,j,k} \Omega^i \otimes \Gamma (S^j E) \otimes \Gamma (\Lambda^k E)
\end{align}
is the connection Lie algebra $\mathfrak{s} \mathfrak{X}_{\nabla}$. The subspace corresponding to $j=k=0$ in the second grouping is precisely the exterior algebra of differential forms $\Omega$. Noting that the actions of $\mathfrak{s}\mathfrak{X}_{\nabla}$ and $\mathfrak{s}\mathfrak{X}$ coincide when restricted to $\Omega$, then altogether we can see that the souped-up Lie algebra $\mathfrak{s}(\Omega \otimes \mathfrak{s}\mathfrak{X}_d, \Omega)$ is a summand of $\mathfrak{s}$. We will refer to this containment when we construct a sheaf of vertex algebras from this souped-up Lie algebra.

\begin{rmk}
The particular element $\mathring{\nabla}_X \triangleq [\nabla, \iota_X] \equiv \nabla\iota_X + \iota_X\nabla$ acts on $\omega \otimes s$ by the beautiful formula
\[ \mathring{\nabla}_X (\omega\otimes s) = \mathcal{L}_X \omega \otimes s + \omega \otimes \nabla s.\]
This is a symmetrized version of the usual covariant derivative $\nabla_X \triangleq \iota_X\nabla$, and can be though of as the covariant analog of the Cartan formula $\mathcal{L}_X = d\iota_X + \iota_X d$. From this formula one can then derive the nice bracket
\[ [\mathring{\nabla}_X, \iota_Y] = \iota_{[X,Y]}.\]
\end{rmk}
\end{example}

\subsection{Souped-up Lie algebra modules}
Upon generating a vertex algebra from a souped-up Lie algebra $\mathfrak{s}$, we will show that $\mathfrak{s}$ has a natural action on the vertex algebra. This will be an instance of a \textit{souped-up Lie algebra module}.

\begin{defn}[Souped-up Lie algebra module]
\label{def:smodule}
A vector space $\mathfrak{v}$ is a module for a souped-up Lie algebra $\mathfrak{s}(\mathfrak{g}, \mathfrak{a})$ if it is a left module for the underlying semidirect sum Lie algebra $\mathfrak{g} \oplus \mathfrak{a}$ and a left module for the algebra $\mathfrak{a}$. That is, we require the two representation $\rho$ and $\sigma$:
\begin{align*}
\rho([s,t], v) &\mapsto \rho(s, \rho(t,v)) - \rho(t, \rho(s,v)) \\
\sigma(ab,v) &\mapsto \sigma (a,\sigma(b,v))
\end{align*}
for $a,b\in \mathfrak{a}$, $s,t \in \mathfrak{s}$, and $v \in \mathfrak{v}$.

In this case, we say $\mathfrak{v}$ is an $\mathfrak{s}$-module.
\end{defn}

Although in $\mathfrak{s}$ we have the relation $[ag,h] = a[g,h] + [a,h]g$, an $\mathfrak{s}$-module does not include the relation $\rho(ag, v) = \sigma(a, \rho(g,v)) + \sigma(\rho(a,v), g)$, because in the final term, $\rho(a,v)$ is not an element of $\mathfrak{a}$, so the notation makes no sense.

\section{Vertex algebras generated by souped-up Lie algebras}
\label{sec:soupedupVA}

We now generate a vertex algebra $V(\mathfrak{s})$ from a souped-up Lie algebra $\mathfrak{s} \equiv \mathfrak{s}(\mathfrak{g}, \mathfrak{a})$. The procedure is first to generate the infinite free algebra $F(\mathfrak{s})$, and then to take the quotient by a particular ideal containing not only the vertex algebra ideal from Definition \ref{def:VA}, but also a set of generators that capture the souped-up structure of $\mathfrak{s}$.

\subsection{Definition of a vertex algebra generated by a souped-up Lie algebra}

\begin{defn}[Souped-up vertex algebra ideal $I(\mathfrak{s})$ and vertex algebra $V(\mathfrak{s})$]
\label{def:SVA}
Let $\mathfrak{s} \equiv \mathfrak{s}(\mathfrak{g},\mathfrak{a})$ be a souped-up Lie algebra and let $F(\mathfrak{s})$ be the infinite free algebra generated by the vector space underlying $\mathfrak{s}$, with grading induced by $|\mathfrak{g}| = 1$ and $|\mathfrak{a}| = 0$. Choosing some function $N(s,t) \geq 0$ on $\mathfrak{s} \times \mathfrak{s}$, the \text{souped-up vertex algebra ideal} $I(\mathfrak{s})$ is generated by the generators from Definition \ref{def:VA} (using $N(s,t)$) and the sets
\begin{description}
\item[Lie algebra $\mathfrak{s}$:] $\textbf{s}\sembrack{s,t} \triangleq s_0t - [s,t]$
\item[algebra $\mathfrak{a}$:] $\textbf{a}\sembrack{a,s} \triangleq a_{-1}s - as$
\item[$\mathfrak{a}$-module:] $\textbf{am}\sembrack{a,b,x} \triangleq (ab)_{-1}x - a_{-1}(b_{-1}x)$
\end{description}
for $a,b\in \mathfrak{a}$, $s,t \in \mathfrak{s}$, and monomials $x \in F(\mathfrak{s})$.

The vertex algebra $V(\mathfrak{s})$ is defined as the quotient $F(\mathfrak{s})/J$ where $J$ is any proper ideal that contains $I(\mathfrak{s})$ and (as is true for $I(\mathfrak{s})$) is decomposable into $|\cdot|$-homogeneous subspaces.
\end{defn}

The requirement that ideal $J \supset I(\mathfrak{s})$ be decomposable into $|\cdot|$-homogeneous subspaces implies that $F(\mathfrak{s}) / J$ retains a $|\cdot|$-grading.

As we mentioned earlier, any structure beyond the underlying unital vector space of $\mathfrak{s}$ is forgotten in $F(\mathfrak{s})$. In contrast, the ideals $I(\mathfrak{v})$ and $I(\mathfrak{s})$ are distinct, the former being properly contained in the latter.

\begin{note}
Regarding the function $N$, we will use $s,t$ to denote elements of $\mathfrak{s}$. Thus $N(s,t)$ is defined only on $\mathfrak{s} \times \mathfrak{s}$.
\end{note}

The elements $\textbf{s}\sembrack{s,t}$ in $I(\mathfrak{s})$ will ensure that the bracket between elements of $\mathfrak{s}$ carries over to $V(\mathfrak{s})$ as the product $\circ_0$. Similarly, the elements $\textbf{a}\sembrack{a,x}$ will transfer the associative commutative product on $\mathfrak{s}$ to the product $\circ_{-1}$. The elements $\textbf{am}\sembrack{a,b,x}$ will ensure that $V(\mathfrak{s})$ is a left module for $\mathfrak{a}$, as we will prove in Theorem \ref{thm:smodule}.

For the record, there is some innocuous redundancy between $\textbf{a}\sembrack{a,s}$ and $\textbf{am}\sembrack{a,b,x}$. For example, one can check that
\[ \textbf{am}\sembrack{a,b,s} =  \textbf{a}\sembrack{ab,s} - \textbf{a}\sembrack{a,bs} - a_{-1}\textbf{a}\sembrack{b,s}  \]
for $a,b \in \mathfrak{a}$ and $s \in \mathfrak{s}$. There is also marginal redundancy between $\textbf{i}\sembrack{x;n}$, $\textbf{s}\sembrack{s,t}$, and $\textbf{a}\sembrack{a,s}$, since $\textbf{i}\sembrack{s;0} = \textbf{s}\sembrack{\1,s}$ and $\textbf{i}\sembrack{s;-1} = \textbf{a}\sembrack{\1,s}$. Lastly, it interesting to note that $\textbf{i}\sembrack{x;n}$ may be replaced by the requirement that $N(s,\1) = 1$ for $s \in \mathfrak{s}$ in $\textbf{c}\sembrack{s,t;n \geq N(s,t)}$. We will not prove or use this fact. The proof uses induction on the length of an element.

\begin{rmk}
\label{rmk:CDRcomparison}
The inclusion of the generators $\textbf{s}\sembrack{s,t}$ in $I(\mathfrak{s})$ likens this construction to the current algebra construction (see \cite{kac-1997}), in which a Lie algebra is used to generate a vertex algebra with the equivalence $[s,t] = s_0t$. This construction is also used to define the chiral de Rham sheaf, in which each vertex algebra captures the Lie superalgebra $\mathfrak{s} \mathfrak{X}$ of vector fields in a semidirect sum with its module $\Omega$. This is akin to our souped-up Lie algebra $\mathfrak{s}(\Omega \otimes \mathfrak{s} \mathfrak{X}, \Omega)$, but without the enhanced ring of coefficients for $\mathfrak{s} \mathfrak{X}$.

Our inclusion of the generators $\textbf{a}\sembrack{a,s}$ and $\textbf{am}\sembrack{a,b,x}$ in $I(\mathfrak{s})$ is where a our construction departs from the construction of the chiral de Rham sheaf. In fact, this set of generators shrinks the vertex algebra in comparison, and in fact \textit{kills} it outright unless, as we discuss in Proposition \ref{prop:Nconstraint}, we compensate by removing generators from $I(\mathfrak{s})$ elsewhere. Specifically, we permit $N(s,t)$ to be greater than 1, whereas the chiral de Rham sheaf uses $N \equiv 1$.

To summarize, in comparison to the construction of the chiral de Rham sheaf, we include the extra generators $\textbf{a}\sembrack{a,s}$ and $\textbf{am}\sembrack{a,b,x}$, but simultaneously toss out some generators $s_nt$ for as many $n$ as we need or want.
\end{rmk}

\begin{thm}
\label{thm:smodule}
The vertex algebra $V(\mathfrak{s})$ is an $\mathfrak{s}$-module (Definition \ref{def:smodule}), with the action given by
\begin{align*}
(s,x) &\mapsto s_0x \\
(a,x) &\mapsto a_{-1}x
\end{align*}
for $a \in \mathfrak{a}$, $s \in \mathfrak{s}$, and $x\in V(\mathfrak{s})$.
\end{thm}

\begin{proof}
We must verify the two equalities
\begin{align*}
[s,t]_0x &= s_0(t_0x) - t_0(s_0x) \\
(ab)_{-1}x &= a_{-1}(b_{-1}x)
\end{align*}
in $V(\mathfrak{s})$.

The second follows directly from $\textbf{am}\sembrack{a,b,x} \in I(\mathfrak{s})$, and the first follows from
\[ \textbf{qa}\sembrack{s,t,x;0,0} - \textbf{s}\sembrack{s,t}_0x \equiv [(s_0t)_0x - s_0(t_0x) + t_0(s_0x)] - (s_0t - [s,t])_0x \]
in $I(\mathfrak{s})$.
\end{proof}

\subsection{Injectivity of \texorpdfstring{$\mathfrak{s}$ into $V(\mathfrak{s}) = F(\mathfrak{s}) / I(\mathfrak{s})$}{\textit{s} into \textit{V(s) = F(s)/I(s)}}}

The feature we most desire in a vertex algebra $V(\mathfrak{s})$ generated from a souped-up Lie algebra $\mathfrak{s}$ is that it contains $\mathfrak{s}$ as a subalgebra, with the bracket $[\cdot,\cdot]$ played by $\circ_0$ and the scalar multiplication by $\mathfrak{a}$ played by $\circ_{-1}$. We know this to be true in two basic cases: when $\mathfrak{g} = 0$, in which case $\mathfrak{s}(0, \mathfrak{a}) \cong \mathfrak{a}$ (Theorem \ref{thm:Ainjectivity}); and $\mathfrak{s}(\mathfrak{g}, \C)$, where $\C$ is the trivial representation for $\mathfrak{g}$. In this case, $\mathfrak{s}$ is really just a central extension of $\mathfrak{g}$ by the ring of coefficients $\C$ (Theorem \ref{thm:Ginjectivity}). However, for a general souped-up Lie algebra $\mathfrak{s}$, it is only \textit{conjectured} that $\mathfrak{s}$ survives intact in $V(\mathfrak{s})$ (Conjecture \ref{conj:injectivity}). We will discuss compelling reasons for this conjecture below.

First we state and prove the two known cases.

\begin{thm}
\label{thm:Ainjectivity}
The map taking the souped-up Lie algebra $\mathfrak{s} (0, \mathfrak{a}) \cong \mathfrak{a}$ to its image in $V(\mathfrak{s}) = F(\mathfrak{s}) / I(\mathfrak{s})$, with $N(s,t) \equiv 0$, is injective. Equivalently, $I(\mathfrak{s}) \cap \mathfrak{s} = \{0\}$.
\end{thm}

This proof would require very little effort if every generator of $I(\mathfrak{s})$ had homogeneous length at least 2. In that case, the \textit{length} decomposition of $F(\mathfrak{s})$ would descend onto $V(\mathfrak{s})$, and then it would suffice to note that $\mathfrak{s}$ and the ideal $I(\mathfrak{s})$ intersect trivially since the former space in the length 1 summand, while $I(\mathfrak{s})$ is orthogonal to that. Thus nothing in $\mathfrak{s}$ is killed in the quotient $F(\mathfrak{s}) / I(\mathfrak{s})$.

Alas, this is not the case, since the generators $\textbf{a}\sembrack{a,s} \equiv a_{-1}s - as$ do not have homogeneous length. (All other generators do have homogeneous length, including $\textbf{s}\sembrack{s,t} \equiv s_0t - [s,t]$, since $[s,t] = 0$ in this setup.)

\begin{proof}[Proof of Theorem \ref{thm:Ainjectivity}]
We will construct a homomorphism $R$ (of vector spaces, not of infinite free algebras) from $F(\mathfrak{s})$ to itself that fixes every element of $\mathfrak{s}$ and projects $I(\mathfrak{s})$ onto the orthogonal complement of the length 1 subspace of $F(\mathfrak{s})$. Since $\mathfrak{s}$ is the length 1 subspace of $F(\mathfrak{s})$, it follows that the intersection $I(\mathfrak{s}) \cap \mathfrak{s}$ is trivial, since apparently any element in the intersection has length 1 (being in $\mathfrak{s}$) and length other than 1 (being in the image of $R$) simultaneously.

To begin, we define the homomorphism of vector spaces $r$ by the rules
\begin{align*}
r(a) &\triangleq a \\
r(a_0b) &\triangleq [a,b] \equiv 0 \\
r(a_{-1}b) &\triangleq ab \\
r(\1_{-1}x) &\triangleq x \\
r(x_{-1}\1) &\triangleq x \\
r(x_ny) &\triangleq r(x)_n r(y) \text{ otherwise},
\end{align*}
for $a,b \in \mathfrak{s} = \mathfrak{a}$ and $x,y \in F(\mathfrak{s})$, and extending linearly.

We then define the projection $R \triangleq r^{\infty}$, the recursive application of $r$. Since $r$ systematically shortens the syntactic strings, deleting any monomial with the substring ``$a_0b$," removing ``$_{-1}$" from ``$a_{-1}b$," eliminating the substrings ``$_{-1}\1$" and ``$\1_{-1}$" altogether, and leaving the string alone otherwise, we see that any element stabilizes after a finite number of applications of $r$. This shows that $R$ is indeed a projection.

The fact that $R$ fixes $\mathfrak{a}$ follows precisely from the first rule.

It remains to show that $R$ projects $I(\mathfrak{s})$ onto the orthogonal complement of the subspace of elements with length 1. It suffices to show that this is true for the generators of $I(\mathfrak{s})$. Indeed, if a generator is killed by $R$, then so is any multiple of that generator. On the other hand, if the image of a generator has length at least 2, then so does the image of any multiple. Altogether, it follows that the entire ideal $I(\mathfrak{s})$ is also mapped under $R$ to an element orthogonal to the subspace of length 1, and hence orthogonal to $\mathfrak{s}$.

We begin with the sets of generators $\textbf{s}\sembrack{a,b} \equiv s_0t - [a,b]$, $\textbf{a}\sembrack{a,b} \equiv a_{-1}b - ab$, $\textbf{am}\sembrack{a,b,x} \equiv (ab)_{-1}x - a_{-1}(b_{-1}x)$, and $\textbf{c}\sembrack{a,b;n}$. Applying $R$, we have
\begin{align*}
R(\textbf{s}\sembrack{a,b}) &= 0 \\
R(\textbf{a}\sembrack{a,b}) &= 0 \\
R(\textbf{am}\sembrack{a,b,x}) &= \textbf{am}\sembrack{a,b,R(x)} \\
R(\textbf{c}\sembrack{a,b;n}) &= \textbf{c}\sembrack{a,b;n}.
\end{align*}
In the third case, by inspection we can see that the image is actually 0 in those cases where $x \in \mathfrak{s}$. Regardless, in all cases above, we see that the image is orthogonal to the subspace $\mathfrak{s}$, since all elements have length $\neq 1$.

Next, we address the generators $\textbf{d}\sembrack{x,y;n} \equiv (x_ny)_{-2}\1 - (x_{-2}\1)_ny - x_n(y_{-2}\1)$. There are the following cases:
\begin{align*}
R(\textbf{d}\sembrack{a,b;0}) &= 0 \\
R(\textbf{d}\sembrack{a,b;-1}) &= (ab)_{-2}\1 - (a_{-2}\1)_{-1}b - a_{-1}(b_{-2}\1) \\
R(\textbf{d}\sembrack{x,\1;-1}) &= R(x)_{-2}\1 - R(x)_{-2}\1 - R(x)_{-1}(\1_{-2}\1) \\
R(\textbf{d}\sembrack{\1,x;-1}) &= R(y)_{-2}\1 - (\1_{-2}\1)_{-1}R(y) - R(y)_{-2}\1 \\
R(\textbf{d}\sembrack{x,y;n}) &= (R(x)_nR(y))_{-2}\1 - (R(x)_{-2}\1)_nR(y) - R(x)_n(R(y)_{-2}\1).
\end{align*}
In each case, the image is again orthogonal to $\mathfrak{s}$.

The proofs for the remaining generators are of a similar spirit, and left to the reader.

This argument was independent of specification of $N(s,t)$. In this case, we may choose $N(s,t) \equiv 0$, which is as strong as possible on $\mathfrak{s} \times \mathfrak{s}$.
\end{proof}

\begin{cor}
If we let $\mathfrak{v}$ be the vector space underlying $\mathfrak{s}$, then this also proves Theorem \ref{thm:Vinjectivity} since $I(\mathfrak{v}) \subset I(\mathfrak{s})$ and hence $I(\mathfrak{v}) \cap \mathfrak{v} \subset I(\mathfrak{s}) \cap \mathfrak {s} = \{0\}$.
\end{cor}

The other basic case for which we can prove that $\mathfrak{s}$ injects into $V(\mathfrak{s})$ is a complement to the first case. In this second case, the Lie algebra may be robust, while the associative commutative unital algebra is played by $\C$.

\begin{thm}
\label{thm:Ginjectivity}
Let $\C$ be regarded as a trivial module for a Lie algebra $\mathfrak{g}$ with coefficients in $\C$. The map taking the souped-up Lie algebra $\mathfrak{s}(\mathfrak{g}, \C)$ to its image in $V(\mathfrak{s}) = F(\mathfrak{s}) / I(\mathfrak{s})$, with $N(s,t) \equiv 1$ on $\mathfrak{s} \times \mathfrak{s}$, is injective. Equivalently, $I(\mathfrak{s}) \cap \mathfrak{s} = \{0\}$.
\end{thm}

\begin{proof}
For this souped-up Lie algebra, since the algebra $\mathfrak{a}$ is played by $\C$, the generators $\textbf{a}\sembrack{a,s}$ and $\textbf{am}\sembrack{a,b,x}$ are easily seen to be combinations of $\1_{-1}s - s$. In this case, we may toss out those sets of generators. If we then choose $N(s,t) \equiv 1$ so that $I(\mathfrak{s})$ contains all elements $s_{\geq 1}t$, it follows that $I(\mathfrak{s})$ is precisely the ideal of relations in the well-known current algebra construction (\cite{kac-1997}, \cite{LL-2007}) with vanishing bilinear form $B$ on $\mathfrak{g}$. In a current algebra, it is known that the entire centrally extended generating Lie algebra survives intact in the quotient $V(\mathfrak{s})$, and the theorem is proved.
\end{proof}

These two theorems leads us to conjecture that this injectivity holds for an arbitrary souped-up Lie algebra.

\begin{conj}
\label{conj:injectivity}
There exists some constant $N$ such that, for any souped-up Lie algebra $\mathfrak{s}$, we have $I(\mathfrak{s}) \cap \mathfrak{s} = \{0\}$ with $N(s,t) \equiv N$.
\end{conj}

For the remainder of this paper, we will assume Conjecture \ref{conj:injectivity}. The following two propositions are necessary conditions that are implied by this conjecture.

\begin{prop}
\label{prop:leftmult}
In general, even if $\mathfrak{s}$ is regarded as a two-sided module for $\mathfrak{a}$ with the relation $as = sa$, the set of elements $\textbf{a}\sembrack{a,s} \equiv a_{-1}s - as \in I(\mathfrak{s})$ may not be enlarged to include the elements $s_{-1}a - as$.
\end{prop}

There are cases where such an addition would not affect the vertex algebra $V(\mathfrak{s})$ adversely. These include the two examples in Theorems \ref{thm:Ainjectivity} and \ref{thm:Ginjectivity}. But there are important cases in which adding the generators corresponding to a right multiplication by $\mathfrak{a}$ actually collapses $V(\mathfrak{s})$ to 0

\begin{proof}[Proof of Proposition \ref{prop:leftmult}]
It suffices to exhibit a souped-up Lie algebra in which the inclusion of a particular element $s_{-1}a - as$ forces $I(\mathfrak{s})$ to include $\1$. The crux is to find three elements $a \in \mathfrak{a}$, $g,h \in \mathfrak{g}$ such that $[h,[g,a]] = \1$.

We choose $\mathfrak{g}$ to be the Lie algebra of vector fields on the real line ${\R}$ and $\mathfrak{a}$ to be the algebra of differentiable functions on ${\R}$. Then letting $b$ be a coordinate function on ${\R}$, we let $a = b^2/2$ and $g$  be the coordinate vector field $\partial/\partial b$. With these choices, we have defined a souped-up Lie algebra in which $[g,[g,a]] = \1$.

Suppose we include the generators $\{ g_{-1}a - ag\}$ in $I(\mathfrak{s})$. Then we can express $\1$ as a combination of elements of $I(\mathfrak{s})$ by the equation
\begin{align*}
\1 &= [g,[g,a]] \\
&= \textbf{qc}\sembrack{a,g;-1}_1 h - \textbf{a}\sembrack{a,g} + (g_{-1}a - ag) - \textbf{e}\sembrack{g_0a,g;1} + \textbf{s}\sembrack{g,a}_0g + \textbf{s}\sembrack{[g,a],g} \\
& \quad - \sum_{k \geq 2} \frac{(-1)^{-1+k}}{k!} \left( \textbf{e}\sembrack{D^{k-1}(g_{-1+k}a),g;1} + \textbf{e}\sembrack{D^{k-2}(g_{-1+k}a),g;0} \right).
\end{align*}
\end{proof}

\begin{prop}
\label{prop:Nconstraint}
For any pair of elements $a \in \mathfrak{a}$ and $g \in \mathfrak{g}$ satisfying $[g,a] = \1$, the generators $\textbf{c}\sembrack{s,t;n \geq N(s,t)}$ must omit $a_ng$ for at least two positive values of $n$. In particular $N \geq 3$ in Conjecture \ref{conj:injectivity}.
\end{prop}

Permitting elements $a_ng$ to survive for some $n \geq 1$ illustrates a difference between the construction at hand and the current algebra construction.

\begin{proof}
As in Proposition \ref{prop:leftmult}, we work with the souped-up Lie algebra in which $\mathfrak{g}$ is the Lie algebra of vector fields on ${\R}$ and $\mathfrak{a}$ is the algebra of differentiable functions. We let $a$ be a coordinate function on ${\R}$, and let $g$ be the coordinate vector field $\partial/\partial a$. Note that $[g,a] = \1$.

Suppose to the contrary that $\textbf{c}\sembrack{a,g;n} \in I(\mathfrak{s})$ for all $n\geq 1$ except for some single positive value $N$. To begin with, this hypothesis implies that $Da - a_{-N-2}(a_Ng) \in I(\mathfrak{s})$, since
\begin{equation}
\label{eq:Da}
Da - a_{-N-2}(a_Ng) = \frac{1}{2}\textbf{qa}\sembrack{a,a,g;-1,-1} + a_{-2}\textbf{s}\sembrack{a,g} + \sum_{\substack{k\geq 1 \\ k\neq N}} a_{-k-2}\textbf{c}\sembrack{a,g;k},
\end{equation}
every expression on the right being in $I(\mathfrak{s})$.

Additionally, our hypothesis implies that $I(\mathfrak{s})$ contains $D^N(a_Ng)$ since it can be written as a combination of elements of $I(\mathfrak{s})$:
\[ D^N(a_N g) = \left(a_0 g + \1 \right) + \textbf{s}\sembrack{g,a} - \textbf{qc}\sembrack{g,a;0} + \sum_{\substack{k\geq 1 \\ k \neq N}} \frac{(-1)^k}{k!} D^k \textbf{c}\sembrack{a,g;k} . \]

This implies furthermore that $(a_Ng)_m x \in I(\mathfrak{s})$ for all $m \leq -N-1$ and $m \geq 0$, and all $x$. Indeed, one can check the identity
\begin{align}
\begin{split}
\label{eq:aNg}
&\binom{m+N}{N} (a_Ng)_m x \\
= &\left( D^N(a_Ng)\right)_{m+N} x - \sum_{k=0}^{N-1} \binom{m+N}{k} \textbf{e}\sembrack{D^{N-1-k} (a_Ng),x;m+N-k},
\end{split}
\end{align}
which expresses $(a_ng)_mx$ as a combination of elements in $I(\mathfrak{s})$. The restrictions on $m$ follow from demanding that the binomial coefficient on the left does not vanish.

Finally, we are able to express the element $\1$ as a combination of elements of $I(\mathfrak{s})$, as seen in the following equation:
\begin{align*}
\1 &= \textbf{s}\sembrack{a,g} - \textbf{e}\sembrack{a,g;1} + \left(Da - a_{-N-2}(a_Ng)\right)_1g + \textbf{qa}\sembrack{a,a_Ng,g;-N-2,1} \\
& \quad + \sum_{k \geq 0} \binom{-N-2}{k}(-1)^k \left( a_{-N-2-k} ((a_Ng)_{1+k}g) - (-1)^N (a_Ng)_{-N-1-k}(a_kg)\right)
\end{align*}
On the right side, the third term is contained in $I(\mathfrak{s})$ due to equation (\ref{eq:Da}), and the terms in the final summation are (multiples of) elements of the form $(a_Ng)_mx$, which are in $I(\mathfrak{s})$ due to equation (\ref{eq:aNg}).
\end{proof}

\section{A functor from unital vector spaces and souped-up Lie algebras to vertex algebras}
\label{sec:uVecttoVA}

Thus far, we have focused on generating a vertex algebra first from a unital vector space and next from a souped-up Lie algebra. But as we have seen, there are many vertex algebras that can arise from either classical structure, depending both on the choice of ideal $J$ used to form the vertex algebra quotient and on the choice of function $N$.

In order to turn this procedure into a functor, we must single out an appropriate choice of $J$ and $N$ so that the resulting vertex algebra behaves well under morphisms. The first choice is easy: just use the vertex algebra ideal $I(\mathfrak{v})$ (or $I(\mathfrak{s})$). For simplicity, we may choose the function $N$ to be the constant posited in Conjecture \ref{conj:injectivity}. In Proposition \ref{prop:functorV}, we will show that these choices lead to a veritable functor from either classical category to the category of vertex algebras.

We begin with the definitions of the unital vector space category $\textbf{uVect}$, the souped-up Lie algebra category $\textbf{SLA}$, and the vertex algebra category $\textbf{VA}$.

\begin{defn}[The category $\textbf{uVect}$]
An object in $\textbf{uVect}$ is a unital vector space. A morphism $\phi \colon \mathfrak{v} \to \mathfrak{v}'$ is a homomorphism such that $\phi(\1_{\mathfrak{v}}) = \1_{\mathfrak{v}'}$.
\end{defn}

\begin{defn}[The category $\textbf{SLA}$]
\label{def:SLAcategory}
An object in the category $\textbf{SLA}$ is a souped-up Lie algebra. A morphism between two souped-up Lie algebras $\mathfrak{s}(\mathfrak{g},\mathfrak{a})$ and $\mathfrak{s}(\mathfrak{g}', \mathfrak{a}')$ is a linear map $\phi$ sending $\mathfrak{g}$ to $\mathfrak{g}'$ and $\mathfrak{a}$ to $\mathfrak{a}'$ such that $\phi(\1_{\mathfrak{s}}) = \1_{\mathfrak{s}'}$, $\phi[s,t] = [\phi(s),\phi(t)]$, and $\phi(as) = \phi(a)\phi(s)$ for $a \in \mathfrak{a}$ and $s,t \in \mathfrak{s}$.
\end{defn}

\begin{defn}[The category $\textbf{VA}$]
An object in $\textbf{VA}$ is a vertex algebra, and a morphism between two vertex algebras $V,V'$ is a linear map $\Phi$ that takes $\1_V$ to $\1_{V'}$ and satisfies $\Phi(x_ny) = \Phi(x)_n\Phi(y)$.
\end{defn}

We now define the map $\mathcal{V}$ from either $\textbf{uVect}$ or $\textbf{SLA}$ to $\textbf{VA}$, and prove that it is a functor.

\begin{prop}[The functor $\mathcal{V}$]
\label{prop:functorV}

The map $\mathcal{V}$ taking $\mathfrak{v} \in \ob(\textbf{uVect})$ to the vertex algebra $V(\mathfrak{v}) = F(\mathfrak{v}) / I(\mathfrak{v})$ with $N(u,v) \equiv 0$, and taking the morphism $\phi \colon \mathfrak{v} \to \mathfrak{v}'$ to the map $\Phi \colon V(\mathfrak{v}) \to V(\mathfrak{v}')$ defined by
\begin{align*}
\Phi(x) &\triangleq \phi(x) + I(\phi(\mathfrak{v})) \text{ for } x\in \mathfrak{v} \subset F(\mathfrak{v})\\
\Phi(x_ny) &\triangleq \Phi(x)_n\Phi(y) \text{ for } x,y \in F(\mathfrak{v}),
\end{align*}
is a functor.

$\mathcal{V}$, defined analogously on $\textbf{SLA}$, but with the constant $N(s,t) \equiv N$ from Conjecture \ref{conj:injectivity}, is also a functor.
\end{prop}

In the case of $\mathcal{V} \colon \textbf{SLA} \to \textbf{VA}$, there might be other viable choices for $N(s,t)$ other than a constant. At the very least, since $\Phi$ necessarily takes $I(\mathfrak{s})$ to $I(\phi(\mathfrak{s}))$, we must have \[ \Phi(\textbf{c}\sembrack{s,t;n}) = \Phi(s_nt) = \Phi(s)_n\Phi(t) = \textbf{c}\sembrack{\Phi(s),\Phi(t);n}, \]
in which case $N(s,t)$ must satisfy $N(\Phi(s),\Phi(t)) \leq N(s,t)$. This is certainly satisfied when $N$ is the constant from Conjecture \ref{conj:injectivity}, but it is possible that $N(s,t)$ is merely constant on each space $\mathfrak{a} \times \mathfrak{a}$, $\mathfrak{a} \times \mathfrak{g}$, and $\mathfrak{g} \times \mathfrak{g}$. In the best case, by Theorems \ref{thm:Ainjectivity} and \ref{thm:Ginjectivity} and Proposition \ref{prop:Nconstraint}, we might be able to choose $N(s,t) \equiv 0$ on the first space, 3 on the second, and 1 on the third. This is an open question.

\begin{proof}[Proof of Proposition \ref{prop:functorV}]
We will sketch the proof on $\textbf{uVect}$. The proof on $\textbf{SLA}$ is the same.

Let us verify first of all that $\mathcal{V}$ is indeed a map from $\hom(\textbf{uVect})$ to $\hom(\textbf{VA})$. Let $\phi$ be a morphism between unital vector spaces $\mathfrak{v}$ and $\mathfrak{v}'$. We must check that $\mathcal{V}(\phi) \equiv \Phi$ is really a morphism of vertex algebras, meaning that $\Phi$ is a homomorphism of infinite free algebras, and that $\Phi$ takes the ideal $I \subset F(\mathfrak{v})$ gets mapped to the ideal $I\subset F(\mathfrak{v}')$. The first statement is clear from the definition of $\Phi$. For the second statement, it suffices to show that generators get taken to generators. Checking this is straightforward (using the comments above when checking this on the generators $\textbf{c}\sembrack{x,y;n}$) and is left to the reader.

Having established that $\mathcal{V}(\phi) \in \hom(\textbf{VA})$, we may proceed by verifying that $\mathcal{V}$ takes the identity morphism $\id_\textbf{uVect}$ to the identity morphism $\id_\textbf{VA}$, and that $\mathcal{V}(\phi \circ \psi) = \mathcal{V}(\phi) \circ \mathcal{V}(\psi)$. These two statements are straightforward to verify, and left to the reader.
\end{proof}

\section{Sheaves of vertex algebras generated by sheaves of \texorpdfstring{$\Omega^0$-modules}{modules over smooth functions}}
\label{sec:sheaf}

From Chapter \ref{sec:uVecttoVA}, we have a functor $\mathcal{V}$ from both $\textbf{uVect}$ and $\textbf{SLA}$ to $\textbf{VA}$, taking $\mathfrak{v}$ to $F(\mathfrak{v})/I(\mathfrak{v})$, where $N(u,v)$ is constant on $\mathfrak{v} \times \mathfrak{v}$ (respectively, $\mathfrak{s}$ to $F(\mathfrak{s})/I(\mathfrak{s})$, and $N(s,t)$ constant on $\mathfrak{s} \times \mathfrak{s}$). In the current chapter, we will extend these associations to sheaves of $\Omega^0$-modules, where $\Omega^0$ is the structure sheaf of smooth functions on a manifold $M$. Given any sheaf $\mathcal{E}$ of $\Omega^0$-modules, we will construct a sheaf $\mathcal{E}^{va}$ of vertex algebras. If $\mathcal{E}$ is additionally a sheaf of souped-up Lie algebras, then $\mathcal{E}^{va}$ is additionally a sheaf of $\mathcal{E}$-modules, following Theorem \ref{thm:smodule}.

This construction will apply to all of the souped-up Lie algebras from Section \ref{sec:SLAexamples} that pertain to structures on a smooth manifold $M$, since they are sheaves of $\Omega^0$-modules (as well as sheaves of souped-up Lie algebras).

Throughout this chapter, $\mathfrak{s}$ will denote a souped-up Lie algebra that is additionally a $\Omega^0(U)$-module, where $U$ is some open set of a smooth manifold $M$. Every result applies also to the more general case in which $\mathfrak{s}$ is only a unital vector space. Unless we say explicitly otherwise, we will tacitly assume that a definition or result applies to both a unital vector space and a souped-up Lie algebra.

\subsection{The supported vertex algebra ideal \texorpdfstring{$K(\mathfrak{s})$}{\textit{K(s)}}}
\label{sec:supportedidealK}

On $\mathfrak{s}$ we have the notion of \textit{support} of the element of $\mathfrak{s}$. We shall extend the notion of support to the infinite free algebra $F(\mathfrak{s})$, and then define an operator that kills elements that have no support.

\begin{defn}[Common support projection operator $\pi$]
\label{def:commonsupport}
Consider the infinite free algebra $F(\mathfrak{s})$ generated by $\mathfrak{s}$. Extending the definition of support of a section to products in the algebra $F(\mathfrak{s})$ by the recursive formula
\[ \supp(x_ny) \triangleq \overline{\Int(\supp x) \cap \Int(\supp y) }\]
and then extending by linearity, we define the \textit{common support projection operator} $\pi$ on $F(\mathfrak{s})$ by
\begin{align*}
\pi (s) &\triangleq x \text{ for }x \in \mathfrak{s} \\
\pi (x_ny) &\triangleq \begin{cases} 0 & \text{ if } \supp(x_ny) = \emptyset \\ x_ny & \text{ otherwise} \end{cases}
\end{align*}
\end{defn}

The interpretation of the operator $\pi$ is that if the length 1 factors of a monomial $x$ have no common support, then $\pi$ kills $x$. Otherwise, $x$ is untouched by $\pi$. This makes $\pi$ indeed a projection operator.

In the next definition, we use the operator $\pi$ to kill monomials in $V(\mathfrak{s})$ (monomials, in the sense of Definition \ref{def:VAmonomiallength}) in which the constituent factors from $\mathfrak{s}$ have no common support. In an example, suppose $f$ and $g$ are two functions on ${\R}$ with no common support. If we were to multiply them pointwise, their product would vanish. We also want the product $f_ng$ to vanish in a vertex algebra. This will eventually be necessary to turn a sheaf of $\Omega^0$-modules into a sheaf of vertex algebras.

\begin{defn}[Supported vertex algebra ideal $K(\mathfrak{s})$]
\label{def:supportedSVAideal}
The \textit{supported vertex algebra ideal} $K(\mathfrak{s}) \subset F(\mathfrak{s})$ is the ideal generated by $I(\mathfrak{s})$ (with constant $N(s,t) \equiv N$) and the additional generators
\begin{description}
\item[common support:] $\textbf{k}\sembrack{x} \triangleq x - \pi(x)$.
\end{description}
\end{defn}

The ideal $K(\mathfrak{s})$ thus contains all monomials the constituent factors of which have no common support, since $\pi$ vanishes on such products. The impact is that when one forms the quotient $F(\mathfrak{s}) / K(\mathfrak{s})$, two monomials whose corresponding factors agree on their common support are now equal.

Upon adding the generators $\textbf{k}\sembrack{x}$ to $I(\mathfrak{s})$, there is the possibility that the resulting ideal $K(\mathfrak{s})$ no longer has trivial intersection with $\mathfrak{s}$, as is true for $I(\mathfrak{s})$ for particular souped-up Lie algebras or when $\mathfrak{s}$ is only a unital vector space by Theorems \ref{thm:Ainjectivity}, \ref{thm:Ginjectivity}, and \ref{thm:Vinjectivity}; and which is conjectured to be true for a general souped-up Lie algebra $\mathfrak{s}$. We are assured by the next theorem that this situation does not occur.

\begin{thm}
\label{thm:Kinjectivity}
The map taking $\mathfrak{s}$ to its image in $V(\mathfrak{s}) = F(\mathfrak{s}) / I(\mathfrak{s})$ is injective. Equivalently,
$K(\mathfrak{s}) \cap \mathfrak{s} = \{0\}$.
\end{thm}

This requires Conjecture \ref{conj:injectivity} in the case that $\mathfrak{s}$ is an arbitrary souped-up Lie algebra.

\begin{proof}
To begin, it follows from the definition of the supported vertex algebra ideal $K(\mathfrak{s})$ that we may write
\[ K(\mathfrak{s}) = I(\mathfrak{s}) + (\textbf{k}), \]
where $(\textbf{k})$ is the ideal generated by the elements $\textbf{k}\sembrack{x}$.

There are three important facts to note: First, by Definition \ref{def:commonsupport}, $\pi$ acts as the identity on the subspace $\mathfrak{s}$. Second, $I(\mathfrak{s})$ is invariant under $\pi$, which follows from the fact that each generator of $I(\mathfrak{s})$ is a sum of monomials with the same factors, and hence $\pi$ either kills or preserves all terms alike. Third, $(\textbf{k})$ is killed by $\pi$, which follows since $\pi$ kills each generator $\textbf{k}\sembrack{x}$ of $(\textbf{k})$.

Then we have
\begin{align*}
K(\mathfrak{s}) \cap \mathfrak{s} &= \pi(K(\mathfrak{s}) \cap \mathfrak{s}) \\
&\subset \pi (K(\mathfrak{s})) \cap \pi (\mathfrak{s}) \\
&= \pi (I(\mathfrak{s}) + (\textbf{k})) \cap \mathfrak{s} \\
&\subset I(\mathfrak{s}) \cap \mathfrak{s} \\
&= \{0\}.
\end{align*}
The final step follows either from Conjecture \ref{conj:injectivity} or Theorems \ref{thm:Ainjectivity}, \ref{thm:Ginjectivity}, or \ref{thm:Vinjectivity}.
\end{proof}

\subsection{From presheaves to sheaves of vertex algebras}
\label{sec:VAsheaf}

As before, we will continue tacitly to acknowledge that the forthcoming definitions and results apply equally to souped-up Lie algebras and unital vector spaces, unless we state otherwise.

Having established the viability of the supported vertex algebra ideal $K(\mathfrak{s})$, we use it to define the presheaf $\mathcal{E}^{va}$ of vertex algebras associated to a sheaf of $\Omega^0$-modules $\mathcal{E}$. Using $K(\mathfrak{s})$ instead of $I(\mathfrak{s})$ will be of utmost importance when it comes time to converting this presheaf into a sheaf.

\begin{defn}[Presheaf $\mathcal{E}^{va}$ of vertex algebras]
\label{def:VAsheaf}
For a sheaf $\mathcal{E}$ of $\Omega^0$-modules, we define the presheaf $\mathcal{E}^{va}$ of vertex algebras by the following data. On an open set $U$ we define
\[ \mathcal{E}^{va}(U) \triangleq \mathcal{V}(\mathcal{E}(U)) \equiv F(\mathcal{E}(U)) / K(\mathcal{E}(U))\]
where $\mathcal{V}$ is the functor featured in Proposition \ref{prop:functorV}.

A restriction morphism $\res_{U,U'} \colon \mathcal{E}(U) \to \mathcal{E}(U')$ induces a morphism
\[ \Res_{U,U'} \colon F(\mathcal{E}(U)) \to F(\mathcal{E}(U')) \]
by the formula
\begin{align*}
\Res_{U,U'} x &\triangleq \res_{U,U'}x \text{ for } x\in \mathcal{E}(U) \subset F(\mathcal{E}(U))\\
\Res_{U,U'}(x_ny) &\triangleq (\Res_{U,U'}x)_n(\Res_{U,U'}y) \text{ for } x,y \in F(\mathcal{E}(U))
\end{align*}

Since $\Res_{U,U'}$ also takes $K(\mathcal{E}(U))$ to $K(\mathcal{E}(U'))$, it follows that $\Res_{U,U'}$ descends to the quotient $F(\mathcal{E}(U)) / K(\mathcal{E}(U))$. Then a restriction morphism
\[\Res_{U,U'} \colon \mathcal{E}^{va}(U) \to \mathcal{E}^{va}(U') \]
is induced by $\res_{U,U'}$ according to the commutative diagram
\[
\begin{tikzpicture}
\matrix(m)[matrix of math nodes,
row sep=3em, column sep=3em,
text height=1.5ex, text depth=0.25ex]
{ \mathcal{E}(U) & (F(\mathcal{E}(U)), K(\mathcal{E}(U))) & \mathcal{E}^{va}(U)  \\ \mathcal{E}(U') & (F(\mathcal{E}(U')), K(\mathcal{E}(U'))) & \mathcal{E}^{va}(U') \\};
\path[->]
	(m-1-1) edge node[left]{$\res_{U,U'}$} (m-2-1)
	(m-1-2) edge node[left]{$\Res_{U,U'}$} (m-2-2)
	(m-1-3) edge node[left]{$\Res_{U,U'}$} (m-2-3)
	(m-1-1) edge (m-1-2)
	(m-1-2) edge (m-1-3)
	(m-2-1) edge (m-2-2)
	(m-2-2) edge (m-2-3);
\end{tikzpicture}
\]
\end{defn}

We could very well have defined a presheaf of more general vertex algebras (not necessarily $\Omega^0$-modules), ignoring the support generators $\textbf{k}\sembrack{x}$ altogether. However, such a presheaf cannot in general be made into a sheaf, for otherwise the uniqueness axiom is violated, as we will see in the proof of Theorem \ref{thm:VAsheaf}.

We now prove that the presheaf $\mathcal{E}^{va}$ is a sheaf, beginning with a technical lemma.

\begin{lem}
\label{lem:bumpsupport}
Let $U \subset M$ be open and let $T \subset U$ be an open subset such that $\overline{T} \subset U$. Let $\sigma$ be a smooth bump function on $U$ that has height 1 on $T$ and vanishes smoothly at the boundary of $U$. Lastly, consider an element $x = x\sembrack{s, t, \ldots, u} \in F(\mathcal{E}(U))$ with $s,t,\ldots, u \in \mathcal{E}(U) \subset F(\mathcal{E}(U))$. Then using the notation
\[ \sigma \ast x \triangleq x\sembrack{\sigma s, \sigma t, \ldots, \sigma u}, \]
we have
\[ \supp(x - (\sigma \ast x)) \subset U\backslash T .\]
\end{lem}

\begin{proof}
The proof is an induction on the length of $x$. As the base case, we assume $x$ itself has length 1, and therefore is in $\mathcal{E}(U)$. Then
\[ x - (\sigma \ast x) \equiv x - \sigma x  = (\1-\sigma) x. \]
The factor $\1- \sigma$, and therefore the product $(\1-\sigma)x$, clearly has support only in $U\backslash T$.

For the inductive case, suppose that $\supp (x - (\sigma \ast x)) \subset U\backslash T$ for all elements $x$ with length not exceeding $p$. Without loss of generality, consider a monomial element $w$ of length $p+1$. $w$ necessarily factors as $w = x_ny$ for some elements $x$ and $y$ with lengths not exceeding $p$. Then we have
\begin{align*}
x_ny - (\sigma \ast (x_ny)) &= x_ny - (\sigma \ast x)_n (\sigma \ast y) \\
&= \left[ x_ny - x_n (\sigma \ast y) \right] + \left[ x_n(\sigma \ast y) - (\sigma \ast x)_n (\sigma \ast y) \right] \\
&= x_n(y - (\sigma \ast y)) + (x - (\sigma \ast x))_n (\sigma \ast y).
\end{align*}

Then taking the support of both sides, we have
\begin{align*}
\supp \left( x_ny - (\sigma \ast (x_ny))\right) &= \supp \left( x_n(y - \sigma \ast y) + (x - (\sigma \ast x))_n (\sigma \ast y) \right) \\
&= \supp \left( x_n(y - \sigma \ast y)\right) \cup \supp \left((x - (\sigma \ast x))_n (\sigma \ast y) \right) \\
&\subset \left( \supp x \cap (U\backslash T) \right) \cup \left(  (U\backslash T) \cap \supp (\sigma \ast y) \right) \\
&\subset U\backslash T.
\end{align*}
\end{proof}

\begin{cor}
\label{cor:bumpsupport}
Let $U,T$ and $\sigma$ be as in Lemma \ref{lem:bumpsupport}. Let $\rho$ be a function on $U$ with support in $T$. Then in the quotient $F(\mathcal{E}(U)) / K(\mathcal{E}(U))$, we have, for all $x$,
\[ \rho_n x = \rho_n (\sigma \ast x).\]
\end{cor}

\begin{proof}
We see that the support of the difference between the two sides is
\begin{align*}
\supp \left(\rho_nx - \rho_n(\sigma \ast x)\right) &= \supp \left( \rho_n (x - (\sigma \ast x))\right) \\
&\subset \overline{ \Int (\supp \rho) \cap \Int (\supp (x - (\sigma \ast x)))} \\
&= \overline{T \cap (U\backslash T)} \\
&= \emptyset.
\end{align*}
Thus by Definition \ref{def:commonsupport}, $\pi(\rho_nx - \rho_n(\sigma \ast x)) = 0$, and finally
\[ \rho_nx - \rho_n(\sigma \ast x) = \rho_nx - \rho_n(\sigma \ast x) + \pi(\rho_nx - \rho_n(\sigma \ast x)) \in K. \]
\end{proof}

\begin{thm}
\label{thm:VAsheaf}
The presheaf $\mathcal{E}^{va}$ is a sheaf.
\end{thm}

\begin{proof}
We must verify the existence and uniqueness axioms.

$\textbf{Existence.}$ To demonstrate the existence axiom, we must show that given a covering $\{U^i\}$ of any open set $U \subset M$, and given a section $x^i \in \mathcal{E}^{va}(U^i)$ for each $i$ such that the common restrictions agree, meaning
\begin{equation}
\label{eq:commonrestriction}
\Res_{U^i,U^i \cap U^j} x^i = \Res_{U^j,U^i \cap U^j} x^j,
\end{equation}
then there exists a section $x \in \mathcal{E}^{va}(U)$ with $\Res_{U,U^i} x = x^i$ for all $i$.

By the Shrinking Lemma, each open set $U^i$ contains an open set $T^i$ such that $\overline{T^i} \subset U^i$ and $\{T^i\}$ is also an open cover for $U$. Let $\{\rho^i\}$ be a partition of unity subordinate to the open covering $\{T^i\}$, and let $\{\sigma^i\}$ be a set of smooth bump functions (\textit{not} a partition of unity) such that $\sigma^i$ takes the value 1 on $T^i$ and vanishes smoothly at the boundary of $U^i$. That is, starting from the interior of some $T^i$, each function $\rho^i$ dies smoothly as we approach the boundary of $T^i$, while $\sigma^i$ stays constant at 1. Then $\sigma^i$ dies smoothly as we move from the boundary of $T^i$ to the boundary of $U^i$.

The salient feature of the sets $U^i, T^i$ and the functions $\rho^i, \sigma^i$ is that they satisfy the conditions for Corollary \ref{cor:bumpsupport}, which we will use below.

The global section $x \in \mathcal{E}^{va}(U)$ we seek is
\[ x \triangleq \sum_j \rho^j_{-1} \left(\sigma^j \ast x^j\right).\]
Ignoring the ``$\sigma^j \ast$" for a moment, this definition resembles the usual breakdown of a smooth function into a sum of smooth functions supported on the various open sets $U^j$. The reason we must incorporate the smoothing factors $\sigma^j$ is because the elements $\{\rho^j_{-1} x^j\}$ are not contained in $\mathcal{E}^{va}(U)$, whereas each $\rho^j_{-1} \left(\sigma^j\ast x^j\right)$ is. The length 1 factors of the former elements are $\rho^j$ and various elements of $\mathcal{E}(U^j)$ which are not in general directly multiplied by $\rho^j$ (in the multiplication in $\mathcal{E}(U^j)$), and therefore do not vanish smoothly themselves at the boundary of $U^j$.

We check that $x$ restricts to $x^i$ on each open set $U^i$ as it ought:
\begin{align*}
\Res_{U,U^i} x &= \Res_{U,U^i} \sum_j \rho^j_{-1} \left(\sigma^j \ast x^j\right) \\
&= \Res_{U,U^i} \sum_j \rho^j_{-1} x^j \\
&= \Res_{U,U^i} \sum_j \rho^j_{-1} x^i \\
&= \Res_{U,U^i} \1_{-1} x^i \\
&= x^i.
\end{align*}

To go from the first line to the second, we have noted that the $j$th term is supported on $U^i \cap T^j$, enabling us to apply Corollary \ref{cor:bumpsupport}. Going to the third line, we have used the hypothesis that $x^i$ and $x^j$ agree on $U^i \cap T^j$ (equation (\ref{eq:commonrestriction})).

$\textbf{Uniqueness.}$ For uniqueness, we must show that if two elements $x,y \in \mathcal{E}^{va}(U)$ agree when restricted to each $U^i$ in some open covering $\{U^i\}$ of $U$, then $x = y$. By linearity, we may specialize to $y = 0$. The task is then to show that whenever $\Res_{U,U^i} x = 0$ for all $i$, then $x = 0$ as well.

For a fixed index $i$, viewing $\Res_{U,U^i}$ as a map from $F(\mathcal{E}(U))$ to $F(\mathcal{E}(U^i))$, the kernel of $\Res_{U,U^i}$ is spanned by all monomials $z$ such that $\supp(z) \subset (U^i)^c$, the complement of $U^i$ within $U$. Descending to the quotient $\mathcal{E}^{va}(U) = V(\mathcal{E}(U))$, the kernel of $\Res_{U,U^i}$ is any element of the form $z + K \subset F(\mathcal{E}(U))$ with $z$ as before.

Since the inclusion $\supp(z) \subset (U^i)^c$ holds for all $i$, we have
\begin{equation*}
\supp(z) \subset \cap_i (U^i)^c = (\cup_i U^i)^c = U^c = \emptyset.
\end{equation*}
By Definition \ref{def:commonsupport}, this implies that $\pi(z) = 0$ so that $z = z - \pi(z) \in K$, and thus that $z + K = K$, which is 0 in $\mathcal{E}^{va}(U)$. Thus we have shown that if an element is killed by all restriction morphisms $\Res_{U,U^i}$, then that element must be 0.
\end{proof}

\begin{rmk}
In the proof of Uniqueness above, the inclusion $z - \pi(z) \in K$ is precisely where we have relied on the presence of the additional generators $\textbf{k}\sembrack{x} = x - \pi (x)$ in $K$.
\end{rmk}

To conclude this chapter, we note that the sheaf of vertex algebras $\mathcal{E}^{va}$ contains the underlying classical sheaf $\mathcal{E}$ as a subsheaf. This is our analog of the containment of the de Rham complex in the chiral de Rham sheaf.

\begin{thm}
\label{thm:sheafinjectivity}
The sheaf of vertex algebras $\mathcal{E}^{va}$ contains $\mathcal{E}$ as a subsheaf. Moreover, if $\mathcal{E}$ is a sheaf of souped-up Lie algebras, then $\mathcal{E}^{va}$ is a sheaf of $\mathcal{E}$-modules.
\end{thm}

\begin{proof}
The first statement follows from the injectivity of each $\Omega^0$-module $\mathcal{E}(U)$ into $\mathcal{E}^{va}(U)$, in accordance with Theorem \ref{thm:Kinjectivity}. The second statement follows from Theorem \ref{thm:smodule}.
\end{proof}

\section{Chiral vector bundles and chiral differential geometry}
\label{sec:chiralDG}

Finally, we are able to assemble the results of the previous chapters to define a \textit{chiral vector bundle}. We will continue to assume Conjecture \ref{conj:injectivity}.

The later examples of souped-up Lie algebras in Section \ref{sec:SLAexamples} are actually sheaves of modules for the structure sheaf $\Omega^0$ of smooth functions. By the work done in Chapter \ref{sec:sheaf}, each such sheaf $\mathcal{E}$ generates a sheaf of vertex algebras which is moreover a sheaf of $\mathcal{E}$-modules, by Theorem \ref{thm:sheafinjectivity}. We are particularly interested in Example \ref{ex:connectionLA}.

\begin{defn}
Given a $G$-vector bundle with connection $(E,\nabla)$, the \textit{chiral vector bundle} $\mathcal{E}^{ch(E,\nabla)}$ is the sheaf of vertex algebras generated by the sheaf of souped-up Lie algebras
\[ \mathcal{E}^{E,\nabla} \triangleq \mathfrak{s}( \Omega \otimes \Gamma (SE \otimes \Lambda E) \otimes \mathfrak{s} \mathfrak{X}_{\nabla}, \Omega \otimes \Gamma (SE \otimes \Lambda E)). \]

In the notation of the previous chapter, $\mathcal{E}^{ch(E,\nabla)} \equiv (\mathcal{E}^{E,\nabla})^{va}$.
\end{defn}

We use the term \textit{chiral} to highlight its close relation to the chiral de Rham sheaf, as we will see below.

By Theorem \ref{thm:sheafinjectivity}, $\mathcal{E}^{ch(E,\nabla)}$ contains $\mathcal{E}^{E,\nabla}$ as a subsheaf and is a sheaf of $\mathcal{E}^{E,\nabla}$-modules. By the decomposition in equation (\ref{eq:soupconnLAdecomp}), $\mathcal{E}^{E,\nabla}$ itself notably contains the subsheaf
\[ \mathcal{E}^{M\times \C,d} \triangleq \mathfrak{s}( \Omega \otimes \mathfrak{s} \mathfrak{X}_d, \Omega) .\]

In a diagram, we have the sheaf containments
\[
\begin{tikzpicture}
\matrix(m)[matrix of math nodes,
row sep=4em, column sep=4em,
text height=1.5ex, text depth=0.25ex]
{ \mathcal{E}^{ch(M\times \C,d)} & \mathcal{E}^{ch(E,\nabla)} \\ \mathcal{E}^{M\times \C,d} & \mathcal{E}^{E, \nabla} \\};
\path[right hook->]
	(m-2-1) edge (m-1-1)
	(m-1-1) edge (m-1-2)
	(m-2-1) edge (m-2-2)
	(m-2-2) edge (m-1-2);
\end{tikzpicture}
\]
The sheaf $\mathcal{E}^{ch(M\times \C,d)}$ is our analog of the chiral de Rham sheaf. As we have pointed out in Remark \ref{rmk:CDRcomparison}, the vertex algebra relations between the constituent vertex algebras of $\mathcal{E}^{ch(M\times \C,d)}$ and those in the chiral de Rham sheaf do not quite align.

Every feature of the vector bundle $E$ makes an appearance in $\mathcal{E}^{ch(E,\nabla)}$. For example, the connection $\nabla$ on $E$ induces a connection on $\mathcal{E}^{ch(E,\nabla)}(U)$ by
\[ \nabla \colon x \mapsto \nabla_0x. \]
In accordance with Theorem \ref{thm:smodule}, this is a derivation of all products $\circ_n$ in that
\[ \nabla_0(x_ny) = (\nabla_0x)_ny + x_0(\nabla_ny). \]
This map respects the grading decomposition of $\mathcal{E}^{ch(E,\nabla)}$. One can also define a covariant derivative $\nabla_X \equiv \iota_X \nabla$ by
\[ \nabla_X \colon x \mapsto (\iota_X)_0 (\nabla_0 x). \]

The curvature operator $\nabla^2 \equiv \frac{1}{2}[\nabla,\nabla]$ then becomes an operator on $\mathcal{E}^{ch(E\nabla)}(U)$ by any of the following equivalent expressions:
\[ \nabla^2 s = \frac{1}{2} [\nabla,\nabla] s = \frac{1}{2} (\nabla_0\nabla)_0 s = \nabla_0 (\nabla_0 s). \]
When restricted to $\Omega \otimes \Gamma E \subset \mathcal{E}^{E,\nabla}$, we may compute the trace of the curvature and arrive at the usual Chern-Weil map to the cohomology of $M$. An open question is how to extend this to all of $\mathcal{E}^{ch(E,\nabla)}$.

\begin{rmk}
An advantage we enjoy over the chiral de Rham sheaf is that in our analog $\mathcal{E}^{ch(M \times \C, d)}$, the element $d$ gives rise to a global field
\[ d(\zeta) \triangleq \sum_{n\in \Z}\frac{d \circ_n}{\zeta^{n+1}} \]
under the state-field correspondence (see Remark \ref{rmk:fields}). In the chiral de Rham sheaf on the other hand, the corresponding element exists globally only when the underlying manifold is Calabi-Yau \cite{MSV-1999}.
\end{rmk}

The construction of the chiral vector bundle $\mathcal{E}^{ch(E,\nabla)}$ is a starting point for extending \textit{classical} differential geometry to its string theoretic analog.



\addcontentsline{toc}{section}{References}
\nocite{*}
\bibliographystyle{halpha}
\bibliography{thesisbib}
\end{document}